%% file: root.tex
\documentclass[a4,10pt,twocolumn]{scrartcl}
	\usepackage[top=2.5cm, bottom=2.5cm, left=1.6cm, right=1.6cm]{geometry}

\pdfoutput=1

\usepackage[hang]{footmisc}
\setfnsymbol{wiley}

\makeatletter
\renewcommand{\@biblabel}[1]{[#1]\hfill}
\makeatother

\setkomafont{section}{\Large}
\setkomafont{subsection}{\large}
\setkomafont{subsubsection}{\normal}

\usepackage[format=plain,  
labelsep=period,           
font=small,                
labelfont=bf,              
skip=5pt                   
]{caption}
\usepackage{graphicx}          
\usepackage{amsmath,amssymb,amsfonts}
\usepackage{enumitem}
\usepackage{tabularx}
\usepackage{amsthm}
\usepackage{mathptmx}
\DeclareMathAlphabet{\mathcal}{OMS}{cmsy}{m}{n}
\usepackage[colorlinks=true,%
linkcolor=black,%
urlcolor=black,%
citecolor=black]{hyperref}
\usepackage{algorithm}
\usepackage{units}
\usepackage{xfrac}
\usepackage{bbm}
\usepackage{booktabs}
\usepackage{tikz}
\usetikzlibrary{calc}
\usetikzlibrary{fadings}
\usetikzlibrary{arrows, shapes}
\usetikzlibrary{positioning}
\tikzstyle{block} = [draw, rectangle,minimum height=1.5em, minimum width=3em,node distance=0.5cm]
\tikzstyle{vecArrow} = [thick, decoration={markings,mark=at position
	1 with {\arrow[semithick]{open triangle 60}}},
double distance=1.4pt, shorten >= 5.5pt,
preaction = {decorate},
postaction = {draw,line width=1.4pt, white,shorten >= 4.5pt}]
\tikzstyle{coord} = [coordinate, node distance=0.5cm]
\usetikzlibrary{decorations.markings}
\usepackage{authblk}

\newcommand{\R}{\mathbb{R}}
\newcommand{\Rpp}{\mathbb{R}_{+\!\!+}}
\newcommand{\Spp}{\mathbb{S}_{+\!\!+}}
\newcommand{\Z}{\mathbb{Z}}
\newcommand{\N}{\mathbb{N}}
\newcommand{\norm}[1]{\left\Vert #1 \right\Vert}
\newcommand{\deldel}[2]{\frac{\partial #1}{\partial #2}\,}
\DeclareMathOperator*{\argmin}{arg\,min}
\DeclareMathOperator{\conv}{conv}
\DeclareMathOperator{\epi}{epi}
\DeclareMathOperator{\interior}{int}
\newcommand{\D}{\mathcal{D}}
\newcommand{\Dx}{\mathcal{D}_\mathrm{x}}
\newcommand{\Du}{\mathcal{D}_\mathrm{U}}
\newcommand{\Dxj}{\mathcal{D}_\mathrm{xJ}}
\newcommand{\Dj}{\mathcal{D}_\mathrm{J}}
\definecolor{orange}{rgb}{0.850980401039124,0.325490206480026,0.0980392172932625}
\definecolor{blue}{rgb}{0,0.447058826684952,0.74117648601532} 

\theoremstyle{definition}
\newtheorem{defi}{Definition}

\newtheorem{rem}[defi]{Remark}
\newtheorem{assum}[defi]{Assumption}
\theoremstyle{plain}
\newtheorem{thm}[defi]{Theorem}
\newtheorem{cor}[defi]{Corollary}
\newtheorem{lem}[defi]{Lemma}
\theoremstyle{definition}
\newenvironment{pf}{\begin{proof}}{%
	\end{proof}\ignorespacesafterend
}

\newcommand{\myTitle}{Online learning with stability guarantees:\\ A memory-based real-time model predictive controller}

\begin{document}
\title{\myTitle\footnote{This article is an extended version of \cite{SCHWENKEL2018automatica} including all proofs, an application example, and a detailed description of the used algorithm.}}

\date{}
\author[1]{Lukas Schwenkel}    
\author[1]{Meriem Gharbi}
\author[2,3]{Sebastian Trimpe}
\author[1]{Christian Ebenbauer}

\affil[1]{\emph{Institute for Systems Theory and Automatic Control, University of Stuttgart, Germany}}          
\affil[2]{\emph{Institute for Data Science in Mechanical Engineering, RWTH Aachen University, Germany}}
\affil[3]{\emph{Intelligent Control Systems Group, Max
	Planck Institute for Intelligent Systems, Stuttgart, Germany}}        

\maketitle

\let\thefootnote\relax\footnotetext{\emph{Email addresses:} \texttt{schwenkel@ist.uni-stuttgart.de} (Lukas Schwenkel), \texttt{gharbi@ist.uni-stuttgart.de} (Meriem Gharbi), \texttt{trimpe@dsme.rwth-aachen.de} (Sebastian Trimpe), \texttt{ce@ist.uni-stuttgart.de} (Christian Ebenbauer). The authors thank the International Max Planck Research School for Intelligent Systems (IMPRS-IS) for supporting Lukas Schwenkel. This work was supported in part by the Cyber Valley Initiative and the Max Planck Society.}

\begin{abstract}
\input{0_abstract.tex}
\end{abstract}

\input{1_intro.tex}

\input{2_prelim.tex}

\input{3_lin.tex}

\input{4_nonlin.tex}

\input{5_app.tex}

\input{6_con.tex}

\bibliographystyle{abbrv} 
\bibliography{all_references}

\appendix
\input{A_proofs.tex}
\input{A_algo.tex}
\input{A_lipsch.tex}

\end{document}

%% file: 0_abstract.tex
\textbf{Abstract.}
We propose and analyze a real-time model predictive control (MPC) scheme that utilizes stored data to improve its performance by learning the value function online with stability guarantees. 
For linear and nonlinear systems, a learning method is presented that makes use of basic analytic properties of the cost function and is proven to learn the MPC control law and the value function on the limit set of the closed-loop state trajectory. 
The main idea is to generate a smart warm start based on historical data that improves future data points and thus future warm starts. 
We show that these warm starts are asymptotically exact and converge to the solution of the MPC optimization problem.
Thereby, the suboptimality of the applied control input resulting from the real-time requirements vanishes over time.
Numerical show that existing real-time MPC schemes can be improved by storing optimizatising the proposed ng scheme.

%% file: 1_intro.tex
\section{Introduction}
Model predictive control (MPC) is a control strategy that solves at each sampling instant a finite horizon open-loop optimal control problem (see \cite{RAWLINGS2009}). 
As a consequence, at every sampling instant an input sequence and its resulting state trajectory are computed for the whole prediction horizon, thus continuously generating large amounts of optimization data. 
However, only the first portion of the computed input sequence is applied to the system and typically the remaining part is not stored.
This is wasteful and contradicts our human intuition to memorize our decisions when we are facing recurring problems or tasks in order to successively improve them.
Motivated by this, the key question addressed in this article is how to leverage on optimization data by introducing memory and online learning in real-time MPC algorithms in order to improve their performance.

\begin{figure}[t!]
	\centering
	\begin{tikzpicture}[auto,>=latex']
	\node [block, fill=white!85!black, minimum width=2cm] (sys) {\small system};
	\node [block, fill=white!85!black, minimum width=7.5cm, minimum height=1cm, node distance=0.3cm, below=of sys, align=center] (mpc) {\small real-time MPC controller\hspace{0cm}\ \\[1.1cm]};
	\node [block, fill=white!85!black, node distance=0.3cm, below=of mpc, minimum width=2cm] (mem) {\small memory};
	\node [block, fill=white, node distance=0.2cm, minimum width=3.2cm, above of=mpc, xshift=1.6cm, align=center] (sws) {\small temporal warm start};
	\node [block, fill=white, node distance=0.5cm, minimum width=3.2cm, below of=mpc, xshift=1.6cm, align=center] (dbws) {\small spatial warm start};
	\node [block, fill=white, node distance=2.4cm, left of=mpc, align=center, yshift=-0.2cm] (opt) {\small real-time\\\small optimization \\ \small iteration};
	\node[circle, fill=black, inner sep=0.05cm, node distance=0.85cm, left of= mpc, yshift=-0.2cm] (sw3) {};
	\node[circle, fill=black, inner sep=0.05cm, node distance=-0.35cm, right of= mpc, yshift=0.2cm] (sw2) {};
	\node[circle, fill=black, inner sep=0.05cm, node distance=-0.35cm, right of= mpc, yshift=-0.5cm] (sw1) {};
	\draw [->] (dbws.west) |- (sw1) -- (sw3) -- (opt.east);
	\draw (sws.west)  |- (sw2);
	\draw [->] (opt.west) -- +(-0.75,0) |- node[above, pos=0.7] {\small $u(k)$} (sys.west);
	\draw [-] (mpc.east) +(0,-0.1) -- +(0.3,-0.1) |- node [pos=0.7,above] {\small $x(k)$} (sys.east);
	\draw [->] (mpc.east) +(0,-0.1)-- +(-0.25,-0.1) |- (sws.east);
	\draw [->] (mpc.east) +(0,-0.1)-- +(-0.25,-0.1) |- ($(dbws.east)+(0,0.1)$);
	\draw [vecArrow] (mpc.south west) + (0.4,0) |- node[above, yshift=-0.05cm, pos=0.74]{\small optimization}node[below, pos=0.74]{\small data $\D(k)$} (mem.west);
	\draw [<-] (dbws.east)+(0,-0.1) -- +(0.4,-0.1) |- (mem.east);
	\end{tikzpicture}
	\caption{Illustration of the proposed online-learning real-time MPC scheme. The standard real-time MPC controller consists of a (temporal) warm start solution and an optimization iteration. To improve control performance, past optimization data of the MPC are stored in memory. From these, a spatial warm start is constructed, which improves with more data. By taking in each iteration the warm start that results in a lower cost function value, the stability properties of the original real-time MPC are retained.
	} \label{fig:blkdiag}
\end{figure}
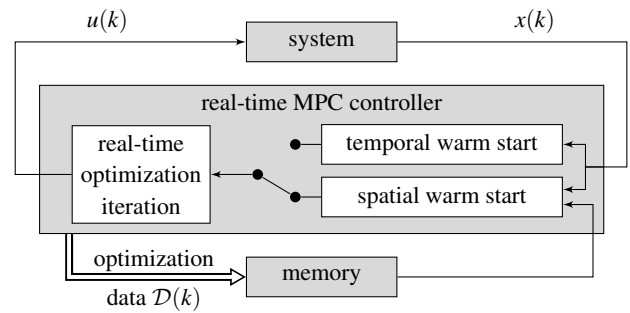

In this work, we will store the previously computed input sequences and utilize them to learn the solution of the MPC optimization problem online, see Fig. \ref{fig:blkdiag}.
The proposed approach specifically targets real-time MPC where, due to a lack of computation time, the MPC optimization problem cannot be solved exactly and thus, suboptimal solutions are applied to the system. 
Hence, there is a need to learn the optimal solution. 

The proposed online learning is based on the idea to store optimization data and leverage it for subsequent optimizations: whenever a new input sequence has to be computed close to a point that was visited before, the optimization iteration is initialized with the previously computed input sequence. 
This way, we can use the available computation time to refine this suboptimal input sequence instead of starting the optimization from scratch. 
At recurrent points, we expect to approach the optimal solution since we improve a suboptimal one over and over again. 
By combining these data-based improvements with the theory of real-time MPC, we will be able to show that online learning is feasible for recurrent points with inherent stability guarantees. 

\textbf{Contributions\ }
The main contributions of this work are as follows. 
We present a real-time MPC scheme based on \cite{FELLER2017} that allows for online learning to improve control performance and that comes with inherent stability guarantees independent of the chosen learning method. 
Further, we introduce novel learning methods tailored to the real-time MPC scheme for both linear and nonlinear systems. 
The learning methods exploit analytic properties of the cost function, namely convexity in the linear case and Lipschitz continuity in the nonlinear case, in order to upper bound the value function. 
This bound is shown to converge to the value function on the $\omega$-limit set of the closed-loop state trajectory. 
Thereby, the input applied to the system converges to the optimal feedback policy. 
The result is an online method that learns the value function and the optimal MPC law asymptotically over the $\omega$-limit set.
Moreover, we present several examples, which illustrate and confirm the results and are used to discuss practical issues and how to circumvent them. 
In particular, we demonstrate how learning and applying the learned control law can be parallelized and executed on different time scales, whereby learning itself does not have to be executed in real time.

\textbf{Related work\ }
This work combines real-time MPC with data-based approaches, both of which being active areas of research. 
The proposed online learning approach builds on top of a real-time MPC scheme and especially its stability properties.
In real-time MPC, there are several results on stability, which all take the suboptimality of the applied input sequence explicitly into account. 
The work \cite{SCOKAERT1999} seeks for a feasible solution without executing any optimization steps, which makes this method unsuitable as basis for learning the optimal solution from its data.
The works \cite{ALAMIR2015}, \cite{BEMPORAD2015}, and \cite{MCGOVERN1999} optimize the input sequence until a certain accuracy is achieved.
Hence, to improve their performance by including learning, we would need to adjust this accuracy according to the learning progress, which is rather difficult to do.
That is why we assume a method for the proposed learning approach that optimizes until a certain computation time is reached like those in \cite{DIEHL2005}, \cite{FELLER2017}, and \cite{ZEILINGER2014}. 
In particular, the real-time MPC scheme presented in \cite{FELLER2017} will be the basis for the online learning method we develop herein. 
In \cite{FELLER2017}, the cost function is established as a Lyapunov function for the coupled system optimizer dynamic, which allows for a straightforward extension of the scheme to learning while preserving stability.

Learning approaches in MPC are becoming increasingly popular. Different ways for incorporating data in MPC have been suggested in literature and can be summed up in two main categories: \emph{(i)} system input-output data to learn the model; and \emph{(ii)} sampling data of the MPC policy to approximate it with an explicit control law.
Model learning \emph{(i)} is an important topic of current research since the accuracy of the model has significant impact on the control performance. 
The model can be learned offline, as for example in \cite{CANALE2014}, \cite{KOCIJAN2004}, \cite{LENZ2015}, \cite{PICHE2000}, and \cite{WILLIAMS2017}, but also online, \cite{ASWANI2013}, \cite{BOUFFARD2012}, \cite{CHOWDHARY2013}, \cite{OSTAFEW2014}, or first learned offline and then refined online \cite{GU2002}, \cite{LIMON2017}. 
Nevertheless, only few of these works, namely \cite{ASWANI2013}, \cite{CANALE2014}, \cite{CHOWDHARY2013}, and \cite{LIMON2017}, establish stability results for their methods. 
The learning of an explicit MPC control law \emph{(ii)} is done offline in order to substitute the online optimization in MPC with an online evaluation of the explicit control law. 
The majority of works in this direction employ neural networks as approximate controllers \cite{AKESSON2006}, \cite{CHEN2018}, \cite{DRGONA2018}, \cite{HERTNECK2018}, \cite{HIROSE2018}, \cite{PARISINI1995}, \cite{PONKUMAR2018} and \cite{ZHANG2016}, while support vector machines are used in \cite{CHAKRABARTY2017}, and the learning problem is formulated as quadratically constrained quadratic program in \cite{DOMAHIDI2011}. 
Out of these works, only \cite{DOMAHIDI2011} and \cite{HERTNECK2018} can guarantee stability.
Since the evaluation of e.g. a neural network is typically much faster than solving the MPC optimization problem, these methods are suited for real-time applications. 
In contrast to our work, these methods learn the control law offline. 
This implies that they have a fixed approximation error and cannot improve online, while our method approaches optimality by learning online. 
By their design, offline methods cannot adapt to changes and need training data beforehand, while our method can be started without prior data.

There are few other works on learning in MPC that do not fit in the two categories.
One of them is \cite{ROSALIA2018}, in which past inputs and the past state trajectory are used for iterative tasks to learn a terminal set and terminal cost of a finite horizon MPC controller such that infinite horizon performance is optimized while stability is ensured.
In \cite{KLAUCO2019}, a neural network is trained offline to identify active constraints in order to warm start an active set algorithm in embedded MPC while guaranteeing online stability.
Although this work also considers learning warm starts, it is conceptually different from ours, since it learns offline and does not improve performance by decreasing the cost.

This article proposes a novel way of leveraging data in MPC, which does not match any of the aforementioned categories and works. Here, we use internal data of the \emph{real-time} MPC algorithm, specifically predicted state and input sequences. 
With this data, we improve the suboptimal real-time controller over time by providing better warm start solutions based on past optimization solutions. 
To the best of the authors' knowledge, this is the first work on learning in real-time MPC that utilizes this internal data.

\textbf{Outline\ }
This article continues with introducing preliminary results, the problem formulation, and the core idea of this work in Section \ref{sec:prelim}. 
Sections \ref{sec:linMPC} and \ref{sec:nonLinMPC} then contain the main results of this article: the online learning method is developed for linear and nonlinear MPC, respectively, and their properties are analyzed in theory and illustrated through numerical examples. 
Section \ref{sec:App} underlines the relevance of the results by an application example. Our conclusions, stated in Section \ref{sec:Conc}, complete the article.

\textbf{Notation\ }
Throughout the article, we will use $\R_+$, $\Rpp$ and $\N_{+}$ to denote the sets of non-negative real, strictly positive real, and strictly positive natural numbers. 
The set of positive definite matrices of dimension $n\in \N_+$ is denoted with $\Spp^n$. For $M\in\Spp$, we define $\norm{x}_M := (x^\top M x)^{1/2}$. 
A continuous function $\alpha : \R_+ \to \R_+$ with $\alpha (0) = 0$ is a $\mathcal K_\infty$-function if  it is strictly increasing and $\alpha (s) \to \infty$ as $s\to \infty$. 
The interior of a set $A\subset \R^n$ is denoted $\interior A$ and the convex hull is $\conv A$. 
For a function $f:\R^n\to \R$ the epigraph is
\begin{align} 
\epi f=\{ (x,\mu)\in\R^n\times \R | f(x)\leq \mu \}.
\end{align}

%% file: 2_prelim.tex
\section{Problem setting: Learning in real-time MPC} \label{sec:prelim}
In this section, we set the stage for developing the main results of this article.  We first introduce the real-time MPC framework \cite{FELLER2017},
in which our learning method will be embedded. We then present the main idea of how to leverage data in this setting, make the learning problem precise, and provide preliminary stability results.
\subsection{Real-time MPC framework} \label{sec:probSetup}
We consider the control of time-invariant discrete-time systems of the form
\begin{align}\label{eq:sysx}
x(k+1)=f\big (x(k),u(k)\big ) + w(k)
\end{align}
where $x(k)\in\R^n$ denotes the state vector, $u(k)\in\R^m$ the control input vector, and $w(k)\in\R^n$ some external signal vector, all at time instant $k\in\N$. We will consider both linear and nonlinear system dynamics $f:\R^n \times \R^m \to \R^n$. In MPC, a typical control task is regulation of the system state to the origin while minimizing some cost function. This is handled by solving at each time instant $k$ an open-loop optimal control problem of the form
\begin{align}
	&\hspace{-2.6cm}J_N\big (U,x(k)\big )=\sum_{j=0}^{N-1} l(x_j,u_j)+F(x_N) \nonumber \\
	J_N^*\big (x(k)\big )=\min_U\ & J_N\big (U,x(k)\big )\label{eq:mpc_prob} \\\text{s.t. }  & x_{j+1}=f(x_j,u_j),\quad j=0,\dots,N-1 \nonumber\\&  x_0=x(k) \nonumber
\end{align}
where $U=[u_0^\top, \dots , u_{N-1}^\top ]^\top\in\R^{Nm}$ denotes the stacked control inputs over the finite prediction horizon $N\in\N_{+}$, and $l:\R^n\times \R^m \to \R_+$ and $F:\R^n\to\R_+$ denote the stage and terminal costs, respectively. 
Since the external signal $w(k)$ is unknown over the prediction horizon, we assume $w(k)=0$ and use the nominal system to predict $x_j$.
We will call $J_N^*$ the \textit{value function} and $J_N$ the \textit{cost function}. 
We assume that state and input constraints are taken into account by a barrier or penalty term in the stage and terminal cost. 

Usually in MPC, it is assumed that the system dynamics and the optimization algorithm evolve on different time scales such that the convergence time of the algorithm can be neglected and an instantaneous solution to \eqref{eq:mpc_prob} is available. In real-time iteration schemes, this often unrealistic assumption is dropped and the optimization algorithm is interpreted as a system with its own dynamics. These are coupled with \eqref{eq:sysx} and are given by
\begin{subequations}\label{eq:sysU}
	\begin{align}
	U(k+1)&=\Phi^{i_{\mathrm T}(k)} \big (U(k),x(k) \big )\\
	u(k+1)&=\Pi_0 U(k+1)
	\end{align} 
\end{subequations}
where $\Phi^{i_{\mathrm T}(k)}:\R^{Nm}\times \R^{n} \to \R^{Nm}$ represents the optimization algorithm and $\Pi_0=[I_m \ 0 \ \dots \ 0]\in\R^{m\times Nm}$ is a projection matrix selecting the input to be actually applied. We divide the optimization algorithm $\Phi^{i_\mathrm{T}(k)}$ into a warm start operator $\Psi_{\mathrm {tw}}:\R^{Nm}\times \R^{n} \to \R^{Nm}$ and an optimization update operator $\Psi_{\mathrm{o}}:\R^{Nm}\times \R^n\to\R^{Nm}$, which is iterated $i_{\mathrm T}(k)\in\N_+$ times
\begin{subequations}\label{eq:phi}
	\begin{align}\label{eq:phi0}
	\Phi^{0}(U,x) &= \Psi_\mathrm{tw} (U,x) \\ \label{eq:psi_o}
	\Phi^{i}(U,x) &= \Psi_\mathrm{o} \big (\Phi^{i-1}(U,x),f(x,\Pi_0U)\big ).  
	\end{align}	
\end{subequations}
In this way, an input for the nominal next state $f(x,\Pi_0U)$ is determined. 

The warm start operator is typically generated by a time shift of the previous input sequence appended with some local control law \cite{FELLER2017}. 
Throughout the article, we will therefore refer to $\Psi_{\mathrm {tw}}$ as \emph{temporal} warm start operator. 
If the closed loop is stable for \emph{any} number of optimization algorithm iterations $i_\mathrm{T}(k)$, the real-time scheme is also called \emph{anytime} MPC, \cite{BEMPORAD2015}, \cite{FELLER2017}. 
A generic result on nominal closed-loop stability of an anytime iteration scheme is proven in \cite{FELLER2017} and stated in the following theorem.
\begin{thm}\label{thm:stability}
	Consider the real-time MPC scheme introduced in \eqref{eq:sysx}--\eqref{eq:phi} with $w(k)=0$ and assume that the stage and terminal costs $l$ and $F$ are positive definite. 
	Further assume the existence of $\underline \alpha, \overline \alpha \in \mathcal K_\infty$ such that 
	\begin{align}\label{eq:classKfncts}
	\underline \alpha (\norm{(U,x)}) \leq J_N(U,x)\leq \overline \alpha (\norm{(U,x)})
	\end{align}
	and assume that $\Psi_\mathrm{tw}$ and $\Psi_\mathrm{o}$ fulfill
	\begin{subequations}\label{eq:ineq}
		\begin{align}\label{eq:ineqw}
		J_N\big (\Psi_\mathrm{tw}(U,x),f(x,\Pi_0 U)\big )-J_N(U,x) &\leq -l(x,\Pi_0U)  \\
		J_N(\Psi_\mathrm{o}(U,x),x)-J_N(U,x) &\leq -\gamma (U,x) \label{eq:ineqo}
		\end{align}
	\end{subequations}
	for all $(U,x)\in\R^{Nm} \times \R^n$, where $\gamma:\R^{Nm}\times \R^n \to \R_+$ is a continuous function with $\gamma(U,x)=0 \Leftrightarrow J_N(U,x)=J_N^*(x)$. 
	Then, for any sequence $\{i_\mathrm{T} (0),i_\mathrm{T}(1),\dots \}$, the origin $(U,x)=(0,0)$ is globally asymptotically stable.
\end{thm}
This theorem is stated in \cite{FELLER2017} for linear systems using a relaxed barrier function formalism to ensure the assumptions. For nonlinear systems, this generic result still applies, however, it is much more challenging to ensure \eqref{eq:classKfncts} and \eqref{eq:ineq} in nonlinear MPC. In \cite{DIEHL2005}, a weaker stability result for a slightly different nonlinear real-time MPC framework is presented, yet without satisfying the assumptions \eqref{eq:classKfncts} and \eqref{eq:ineq}.

\subsection{Problem formulation and main idea}
When the system dynamics and the optimization algorithm operate on similar time scales, there is typically only time for a few optimization iterations until the system requires the next input, i.e. $i_{\mathrm T}(k)$ is small. 
Hence the input might be far from optimality resulting in unsatisfying controller performance. 
To solve this issue, we will leverage the data of \emph{all} the previously computed input sequences, which are usually discarded (except the last input sequence generating a temporal warm start). 
We store all past input sequences, and if the system state arrives close to a point in the state space that has been visited before, we use this previously calculated input sequence to generate a warm start solution. 
The main challenges of this approach are to make `close to a point' mathematically precise and to prove convergence.
Hence, we have to choose which stored input sequence (or combinations of multiple ones) to take as warm start solution at a given location such that the learning scheme asymptotically recovers the optimal policy at recurrent points of the closed loop trajectory.

In more detail, we will denote the warm start generated from the optimization data $\D(k)$ by $\Psi_\mathrm{sw}(\,\cdot\,,\D(k)):\R^n \to\R^{Nm}$ and call it \emph{spatial} to distinguish it from the temporal one. The stored data $\D(k)$ includes: the input sequence $U(j)$ for each time $1\leq j\leq k$, the point in state space at which it was calculated $f(x(j-1),\Pi_0 U(j-1))$, and the cost function value that was achieved $J_\mathrm{N}(U(j),f(x(j-1),\Pi_0 U(j-1))$; that is
\begin{subequations}\label{eq:D}
	\begin{align}\label{eq:Dzj}
	\mathcal D(k)&=\D(k-1)\cup \big \{\big (z_k,J_N(z_k)\big )\big\},\quad \D(0)=\emptyset
	\end{align} with $z_k=(U(k),f(x(k-1),\Pi_0U(k-1)))$.
	For the subsequent analysis, we introduce the notation
	\begin{align}
		\mathcal D_\mathrm{x} &:= \{x\in\R^n\,|\,\exists (U,x,J)\in\mathcal D\} \\
		\mathcal D_\mathrm{xJ} &:= \{(x,J)\in\R^{n} \times \R_+\,|\,\exists (U,x,J)\in\mathcal D\} \\
		\mathcal D_\mathrm{U} &:= \{U\in\R^{Nm}\,|\,\exists (U,x,J)\in\mathcal D\}
	\end{align}
\end{subequations}
and we denote an approximation of $J_N^*(x)$ based on the collected data $\D(k)$ by $J_N^\mathrm{a}(x,\D(k))$.

In a nominal MPC stabilization problem, we cannot expect that the system state repeatedly arrives at the same points in the state space except at the origin, where the optimal input is trivial. In a real-world scenario, however, there are disturbances, periodic operation conditions, set point changes, or reference signals such that more points can be visited several times. Thus, \textit{learning} the optimal control at these points is a meaningful task. We generically model such situations with the signal $w(k)$, which influences the shape of the $\omega$-limit set of the sequence of points for which an input is computed,
\begin{align}\label{eq:omega}
\Omega=\big\{y\in\R^n\big|\exists k_i \to \infty: f\big (x(k_i),\Pi_0U(k_i)\big )\to y\big\}.
\end{align}
A point in $\Omega$ is called \emph{limit} or \emph{recurrent} point and is reached infinitely often arbitrarily closely. Only at such points can we expect learning the optimal policy to be possible because the optimization iteration usually converges asymptotically to the optimum and thus needs an infinite number of iterations in general.

With this, we can now make the objective of this work precise.  We aim to design the spatial warm start operator $\Psi_{\mathrm{sw}}$ such that it will converge to the optimal policy for all $x\in\Omega$ as $k\to\infty$, i.e.
\begin{align}\label{eq:goal}
\forall x \in \Omega:\lim_{k\to \infty} J_N\big (\Psi_\mathrm{sw}\big (x,\D(k)\big ),x\big )= J_N^*(x).
\end{align}
Hence, we want to learn the value function $J_N^*$ on $\Omega$ and the corresponding optimal control input.

While learning the optimal control, we do not want to jeopardize stability. Hence, we need to make sure that the warm start solution that is used to initialize the optimization iteration satisfies \eqref{eq:ineqw}. We can achieve stability by exploiting the temporal warm start (which is shown to be stabilizing) and by only applying the spatial warm start when it yields a lower cost. That is, we replace \eqref{eq:phi0} with the new warm start operator
\begin{subequations}\label{eq:psi_w}
\begin{flalign}
&\Phi^0(U,x,\mathcal D)=\begin{cases}
\Psi_\mathrm{tw} (U,x ) & J_{N\mathrm t}<J_{N\mathrm s}\\
\Psi_\mathrm{sw}\big (f(x,\Pi_0U),\mathcal D\big ) & \mathrm{else,}
\end{cases}\hspace{-0.5cm}&\\
&\text{with}\ J_{N\mathrm t}=J_N\big (\Psi_\mathrm{tw}(U,x ), f(x,\Pi_0U )\big )&\\&\phantom{\text{with}\ }J_{N\mathrm s}= J_N\big (\Psi_\mathrm{sw}\big (f(x,\Pi_0U),\mathcal D\big ), f(x,\Pi_0U )\big).\hspace{-1cm}&
\end{flalign}
\end{subequations}
Through this intuitive approach, the resulting online learning scheme is inherently stable as stated in the following corollary.
\begin{cor}
	Consider the real-time MPC scheme introduced in \eqref{eq:sysx}--\eqref{eq:phi}, \eqref{eq:D}, \eqref{eq:psi_w} with $w(k)=0$ and assume \eqref{eq:classKfncts}, \eqref{eq:ineq}. Then the origin $(U,x)=(0,0)$ is globally asymptotically stable for all spatial warm start operators $\Psi_{\mathrm{sw}}$.
\end{cor}
\begin{rem}
	The stability result considers the nominal system with $w(k)=0$ since stability of the origin can not be achieved with a non-vanishing $w(k)$.
	The learning problem, however, is trivial for $\Omega=\{0\}$, thus, we require $w(k_i)\neq 0$ for some $k_i \to \infty$ to render $\Omega\neq\{0\}$.
	Nevertheless, nominal stability is meaningful and essential for the learning task.
	\begin{itemize}
		\item As a first example, consider a scenario as typically studied in iterative learning control, where, after some finite time period, the system is reset to a new initial state. 
		Let the set of initial states be $\mathbb{X}_0$. 
		Then, if we model $w(k_i)=-f(x_{k_i}, u_{k_i}) + \bar x$, $\bar x \in \mathbb X_0$, for a sequence of reset times $\{k_i\}\to\infty$ and $w(k)=0$ for $k\neq k_i$, this would correspond to an iterative learning scenario for multiple initial states in $\mathbb{X}_0$, where the goal is to steer $\bar x\in \mathbb{X}_0$ close to the origin.
		While executing an iteration it is $w(k)=0$, hence, global asymptotic stability of the nominal closed loop system is essential to guarantee that the controller does the right thing (not growing unbounded and approaching the origin).
		This iterative learning interpretation is also discussed in the unicycle example in Section 4 of the paper.
		\item As a second example, let $w(k)$ be Gaussian noise. 
		Then $w(k)$ would render the set of recurrent points nonempty (presumably the whole state space). 
		Nevertheless, the stability result still makes sense in terms of expectations. 
		Since a complex stochastic analysis is beyond the scope of this work, we do not explicitly consider this case. 
		Still, we can think of $w(k)$ as a deterministic noise signal or a realization of a stochastic process. 
		In general, feedback stabilization is essentially only meaningful in the presence of disturbances. 
		Often, in nominal MPC, disturbances are associated to nonzero initial states.
		\item As a third example, consider a tracking of references or set points. 
		These are often based on a nominal stabilizing controller, which is used to stabilize different set points. 
		When the set point changes, this can be interpreted as a new initial condition and modeled by the signal $w(k)$ as in the first example. 
		This scenario is treated in an Application example in Section \ref{sec:App}.
	\end{itemize}
	In general, we assume that $w(k)$ is of such nature that the closed-loop trajectory $(x(k),U(k))$ stays bounded.
	This ensures that $\Omega$ is nonempty due to the Bolzano-Weierstrass theorem and that $\Omega$ is approached by $f (x(k),\Pi_0U(k) )$ as $k\to\infty$.
\end{rem}
Overall, the proposed MPC iteration scheme measures at every time instant the current system state, generates the warm starts, chooses the better one, executes the optimization iteration, and stores its data before it finally applies the first portion of the input sequence to the system. 
A pseudo-code description of this procedure is given in Algorithm \ref{algo:anyLearn}. 
The main focus of this article is line $4$ and $10$ of the algorithm. 
More specifically, we will design in Section \ref{sec:linMPC} and \ref{sec:nonLinMPC} spatial warm start operators for linear and nonlinear MPC, respectively. 
Throughout the article, we assume that a temporal warm start satisfying the conditions in Theorem \ref{thm:stability} is given.
\begin{algorithm}
	\caption{Real-time MPC scheme with learning}\label{algo:anyLearn}
    \vspace{-0.2cm}
	\begin{enumerate}[label={\small \arabic*:}]\setlength{\itemsep}{-1ex}
		\item \textbf{for} $k=1,2,\dots $ \textbf{do}
		\item \quad \emph{obtain current state:} $x(k)$
		\item[ ]\quad \emph{generate warm starts:}
		\item \quad $U_\mathrm{t}=\Psi_\mathrm{tw}(U(k),x(k))$ \emph{(temporal)} 
		\item \quad $U_\mathrm{s}=\Psi_\mathrm{sw}(f(x(k),\Pi_0U(k)),\D(k))$ \emph{(spatial)}
		\item \quad $U_+=U_\mathrm{t}$ \textbf{if} $J_{N\mathrm t}<J_{N\mathrm s}$ as per \eqref{eq:psi_w} \textbf{else} $U_+=U_\mathrm{s}$
		\item[ ]\quad \emph{real-time optimization iteration:}
		\item \quad \textbf{for} $i=1,2,\dots,i_{\mathrm T}(k)$ \textbf{do} \emph{optimizer update}
		\item \qquad $U_+= \Psi_{\mathrm{o}}(U_+,f(x(k),\Pi_0U(k)))$
		\item \quad \textbf{end for}
		\item \quad $U(k+1)=U_+$
		\item[ ]\quad \emph{memorize data:}
		\item \quad $\mathcal D(k+1)=\D(k) \cup \{(z_k,J_N(z_k))\}$ as per \eqref{eq:D}
		\item \quad \emph{apply first input:} $u(k+1)=\Pi_0 U(k+1)$
		\item \textbf{end for}
	\end{enumerate}\vspace{-0.2cm}
\end{algorithm}

%% file: 3_lin.tex
\section{Leveraging data in real-time linear MPC}\label{sec:linMPC}
In this section, we assume \eqref{eq:sysx} to be linear
\begin{align}\label{eq:sysx_lin}
f(x,u)=Ax+Bu,
\end{align}
with $A\in\R^{n\times n},\ B \in \R^{n\times m}$, the stage and terminal cost are quadratic, and the polytopic state and input constraints are incorporated in the costs via relaxed logarithmic barrier functions $\hat B_\mathrm{x}:\R^n\to\R_+$, $\hat B_\mathrm{u}:\R^m\to\R_+$ (see \cite{FELLER2015} or \cite{FELLER2017b} for definition). Under these assumptions, $l$ and $F$ in \eqref{eq:mpc_prob} become
\begin{align}\label{eq:l}
l(x,u)&=\norm{x}_Q^2+\norm{u}_R^2+\varepsilon \hat B_\mathrm{x}(x)+\varepsilon \hat B_\mathrm{u} (u),\\
F(x)&=\norm{x}_P^2, \label{eq:F}
\end{align}
with the design parameters $\varepsilon\in\Rpp$, $Q\in\Spp^n$ and $R\in\Spp^m$, as well as $P\in\Spp^n$ resulting from $\varepsilon$, $Q$, $R$ and the constraints (see \cite{FELLER2017b} for details). 
In \cite{FELLER2017}, a temporal warm start and an optimization update operator are defined such that the conditions \eqref{eq:classKfncts} and \eqref{eq:ineq} of Theorem \ref{thm:stability} are satisfied. 
Moreover, it has been shown in \cite{FELLER2017} that constraint satisfaction can be guaranteed with a finite barrier parameter $\varepsilon$ if the relaxation of the logarithmic barrier functions and the initialization of the scheme are chosen suitably. 
We will refer to this scheme as linear anytime MPC. 
For this setting, we are going to present a method to learn the optimal policy and the value function in the sense of \eqref{eq:goal}, analyze its convergence properties, and demonstrate the method in an example. 
\subsection{A spatial warm start based on convexity}
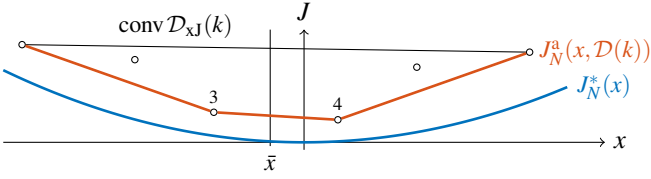
\begin{figure}
	\centering
	\resizebox{\linewidth}{!}{
		\begin{tikzpicture}
		\draw[->] (-4,0)--(4,0) node [right] {$x$};
		\draw[->] (0,-0.1)--(0,1.5) node [above] {$J$};
		\node [inner sep=0.03cm, circle, draw] at (1.5,1) {};
		\node [inner sep=0.03cm, circle, draw] (e4) at (3,1.2) {};
		\node [inner sep=0.03cm, circle, draw] (e2) at (-1.2,0.4) {};
		\node [inner sep=0.03cm, circle, draw] (e1) at (0.45,0.3) {};
		\node [inner sep=0.03cm, circle, draw] (e3) at (-3.75,1.3) {};
		\node [inner sep=0.03cm, circle, draw] at (-2.25,1.1) {};
		\draw (-0.45,-0.05)node [below] {\small $\bar x$} -- (-0.45,1.5);
		\draw (e1)node [above] {\scriptsize $4$} -- (e2)node [above] {\scriptsize $3$} -- (e3) -- node[above,pos=0.3] {\small $\conv \Dxj(k)$} (e4) -- (e1);
		\draw [line width =1pt, orange] (e4) node[right]{\small $J_N^\mathrm{a}(x,\D(k))$} -- (e1) -- (e2) -- (e3);
		\draw [smooth, domain=-4:3.5, samples=40,variable=\x,blue,line width=1pt] plot ({\x},{0.06*\x*\x}) node [right]{\small $J_N^*(x)$};
		\end{tikzpicture}}
	\caption{Illustration of the spatial warm start \eqref{eq:psi_sw_lin} in the $(x,J)$-domain with data points $\Dxj(k)$ (\raisebox{0.03cm}{\scriptsize$\circ$}). The lower boundary (in orange) of their convex hull is an upper bound $J_N^\mathrm{a}(x,\D(k))$ of the value function $J_N^*(x)$ (blue), since $J_N^*$ is convex and optimal. The spatial warm start solution at $\bar x$ is a convex combination of the inputs from points 3 and 4.}\label{fig:psi_sw_lin}
\end{figure}

The main idea for the spatial warm start generation is to exploit convexity of the cost function and use convex combinations of past data points. More formally, we define the spatial warm start at $\bar x\in\conv \Dx(k)$ by
\begin{align}\label{eq:psi_sw_lin}
\big (\Psi_\mathrm{sw}(\bar x,\mathcal D),\bar x,J_N^\mathrm{a}(\bar x,\mathcal D)\big )&=\argmin_{\substack{(U,x,J)\in \conv \mathcal D\\ x=\bar x}} J.
\end{align}
Hereby, we dropped the time dependency of the data $\D=\D(k)$ from \eqref{eq:D} for the sake of clarity. 
Furthermore, $\bar x$ denotes the point at which \eqref{eq:mpc_prob} is to be solved and $J_N^\mathrm{a}(\cdot,\D(k))$ denotes a convex approximation of the value function $J_N^*(x)$ based on the data $\D(k)$. 
A graphical illustration of \eqref{eq:psi_sw_lin} is given in Fig. \ref{fig:psi_sw_lin}. 
Notice that \eqref{eq:psi_sw_lin} is only feasible for $\bar x\in\conv\Dx(k)$ and thus the domain of $J_N^\mathrm{a}(\cdot,\D(k))$ is $\conv\Dx$. 
For $\bar x\notin \conv \Dx(k)$, we cannot compute a spatial warm start in this fashion and have to take the temporal one. 

Geometrically, $J_N^\mathrm{a}$ is a piecewise affine function and partitions $\conv \Dx(k)$ into $n$-simplices, on which it is affine. We will refer to this partition as triangulation of $\conv \Dx(k)$ inspired by topology. For computing the spatial warm start at $x$, one has to find the $n$-simplex that contains $x$, construct the convex combination of its extreme points that leads to $x$, and combine the inputs corresponding to the extreme points in the same way.

In the remainder of this section we will show that this spatial warm start converges to the optimal policy at recurrent points \eqref{eq:omega} by showing that $J_N^\mathrm{a}(\cdot,\D(k))$ upper bounds $J_N(\Psi_{\mathrm{sw}}(x,\D(k)),x)$ if $J_N$ is convex (Lemma \ref{lem:JNaProps}), that $J_N$ is convex (Lemma \ref{lem:JNconv}), and that $J_N^\mathrm{a}(\cdot,\D(k))$ converge to the value function (Theorem \ref{thm:convergence}).

\begin{lem}\label{lem:JNaProps}
	Let $J_N:\R^{Nm}\times \R^n\to\R_+$ be convex with respect to $(U,x)$. Further, let $z_j\in\R^{Nm}\times \R^n,\ j\in\N$ be any sequence, and let $\mathcal D(j)\subset \R^{Nm}\times \R^n\times \R_+$ be the corresponding sequence of data sets defined in \eqref{eq:Dzj}. 
	Moreover, let $\Psi_{\mathrm{sw}}(\cdot,\D(j)):\conv\Dx(j)\to\R^{Nm}$ and $J_N^\mathrm{a}(\cdot,\D(j)):\conv\Dx(j)\to\R_+$ be defined in \eqref{eq:psi_sw_lin} and let $k\in\N_+$. 
	Then (i) $J_N(\Psi_{\mathrm{sw}}(x,\mathcal D(k)),x)\leq J_N^{\mathrm{a}}(x,\mathcal D(k))$ and (ii) $J_N^\mathrm{a}(\cdot,\D(k))$ is convex.
\end{lem}
For a proof, see Appendix \ref{sec:appendix_proofs}.
In order to apply Lemma \ref{lem:JNaProps}, convexity of $J_N$ is needed, which is established for the linear anytime MPC by the following lemma.
\begin{lem}\label{lem:JNconv}
	The cost function $J_N:\R^{Nm}\times \R^n \to \R_+$ of linear anytime MPC, i.e. \eqref{eq:mpc_prob} with \eqref{eq:sysx_lin}, \eqref{eq:l}, \eqref{eq:F} is convex in $(U,x)$.
\end{lem}
The proof is given in Appendix \ref{sec:appendix_proofs}.
In view of Lemma \ref{lem:JNaProps} and Lemma \ref{lem:JNconv}, we can show the following convergence theorem, which establishes that we can indeed asymptotically learn the value function and the MPC law in $\Omega$. 
\begin{thm}\label{thm:convergence}
	Consider the system \eqref{eq:sysx_lin} controlled by \eqref{eq:mpc_prob}--\eqref{eq:D}, \eqref{eq:psi_w}, \eqref{eq:l}--\eqref{eq:psi_sw_lin} and $\Omega$ from \eqref{eq:omega}. 
	Then for all $x \in \Omega$ with $x\in \mathrm{int}\, \conv \Dx(k_0)$ for  some $k_0\in\N$ it holds
	\begin{align}\label{eq:convergenceSW}
	\lim_{k\to\infty} J_N \big (\Psi_\mathrm{sw}\big (x,\mathcal D(k) \big ),x\big )= J_N^*(x).
	\end{align}
\end{thm}
\begin{pf} Step 1) $J_N^\mathrm{a}(\cdot,\D(k))$ converges: For $k\geq k_0$, it holds $\conv \D(k_0)\subseteq \conv \D(k)$ and hence $J_N^\mathrm{a}$ is defined on $\conv \D(k_0)$ for all $k\geq k_0$. 
	Further, since $\conv \D(k)\supseteq \conv \D(k-1)$ and the minimum over a larger set can only be smaller, we can conclude that $J_N^\mathrm{a} (y,\mathcal D(k) )\leq J_N^\mathrm{a} (y,\mathcal D(k-1) )$ for all $y\in \conv \Dx(k_0)$, $k\geq k_0$. 
	In order to apply Lemma \ref{lem:JNaProps}, $J_N$ has to be convex, which is shown in Lemma \ref{lem:JNconv}.
	Due to (i) in Lemma \ref{lem:JNaProps} and the optimality of $J_N^*$ we have
	\begin{align}\label{eq:sandwich}
	J_N^*(y)\leq J_N(\Psi_{\mathrm{sw}}(y,\mathcal D(k)),y)\leq J_N^{\mathrm{a}}(y,\mathcal D(k))
	\end{align}
	for all $y\in \conv \Dx(k_0)$, $k\geq k_0$. Hence, $J_N^\mathrm{a} (\cdot,\D(k))$ is nonincreasing and bounded from below on $\conv \Dx(k_0)$ and thus converges pointwise $\lim_{k\to \infty} J_N^\mathrm{a}(\cdot,\D(k))=J_N^\infty(\cdot)$ on $\conv \Dx(k_0)$.
	Due to (ii) in Lemma \ref{lem:JNaProps}, this is a sequence of convex functions $J_N(\cdot,\D(k))$, which has a convex limit $J_N^\infty(\cdot)$. 
	Every convex function is locally Lipschitz continuous on open subsets (see e.g. \cite{BORWEIN2006}), hence $J_N^\infty(\cdot)$ is locally Lipschitz continuous on $\interior \conv \Dx(k_0)$. 

	Step 2) $J_N^\mathrm{a}(x,\D(k))$ converges to $J_N^*(x)$: Let $x\in\Omega\,\cap\,\mathrm{int}\, \conv \Dx(k_0)$ and let $\mathcal C\subset \interior \conv \Dx(k_0)$ with $\interior \mathcal C\ni x$ be compact, then there exists a sequence $k_i \to \infty$ such that $\xi_i=f\big (x(k_i),\Pi_0U(k_i)\big )\to x,\ \xi_i \in \mathcal C.$
	Let $L_k$ be the Lipschitz constant of $J_N^\mathrm{a}(\cdot,\D(k))$ on $\mathcal C$ and $L_\infty$ for $J_N^\infty(\cdot )$, respectively, which are finite since $\mathcal C$ is compact and the functions are locally Lipschitz. 
	Hence, $L_k\to L_\infty$ as $k\to \infty$ is a real valued converging sequence and is therefore upper bounded by some $M\in\R$. It follows
	\begin{align*}
	\big|J_N^\mathrm{a} \big (\xi_i, \D(k_{i}&+1)\big )-J_N^\mathrm{a}\big (\xi_i,\D(k_i)\big ) \big| \leq 2M \norm{\xi_i-x}\\&+ \big|J_N^\mathrm{a} \big (x, \D(k_{i}+1)\big )-J_N^\mathrm{a}\big (x,\D(k_i)\big ) \big| \to 0.
	\end{align*}
	Thus the following chain of inequalities
	\begin{align}\label{eq:sandwich2}
	\begin{split}
	J_N^{\mathrm a} \big (&\xi_i,\D(k_i+1)\big )\leq J_N\big (U(k_i+1),\xi_i\big )\\&\leq J_N\big (\Phi^1\big (U(k_i),x(k_i),\D(k_i)\big ),\xi_i\big )\\&\leq J_N\big (\Phi^0\big (U(k_i),x(k_i),\D(k_i)\big ),\xi_i\big )\\&\leq J_N\big (\Psi_\text{sw}\big (\xi_i,\D(k_i)\big ),\xi_i\big )\  \leq\ J_N^{\mathrm a}\big (\xi_i,\D(k_i)\big )
	\end{split}
	\end{align}
	is a chain of equalities in the limit, where the first inequality holds due to \eqref{eq:D} and \eqref{eq:psi_sw_lin} and the other inequalities due to \eqref{eq:ineqo}, \eqref{eq:ineqo}, \eqref{eq:psi_w} and Lemma \ref{lem:JNaProps} (i) in this specific order. 
	Therefore as $i\to\infty$ the decrease of the optimizer update operator \eqref{eq:ineqo} $\gamma \left( \Phi^0\big (U(k_i),x(k_i),\D(k_i)\big ), \xi_i \right) \to 0$, which is due to $\gamma(U,x)=0 \Leftrightarrow J_N(U,x)=J_N^*(x)$ only possible if 
	\begin{align*}
		J_N \left( \Phi^0\big (U(k_i),x(k_i),\D(k_i)\big ), \xi_i \right)-J_N^*(\xi_i) \to 0.
	\end{align*}
	This is due to \eqref{eq:sandwich2} equivalent to $J_N^\mathrm{a} \big ( \xi_i,\D(k_i) \big )-J_N^*(\xi_i) \to 0.$
	Since $J_N^*$ is convex, it is also locally Lipschitz continuous on $\interior \conv \Dx(k_0)$ and we assume that $M$ is a Lipschitz constant of $J_N^*$ on $\mathcal C$ (if not choose $M$ large enough). 
	Hence, it follows
	\begin{align*}
	\big | J_N^\mathrm{a}& \big ( x,\D(k_i) \big ) -J_N^*(x) \big | \leq \big | J_N^\mathrm{a} \big ( \xi_i ,\D(k_i) \big ) -J_N^*(\xi_i) \big |\\
	&+\big | J_N^\mathrm{a} \big ( x ,\D(k_i) \big )-J_N^\mathrm{a} \big ( \xi_i ,\D(k_i) \big ) + J_N^*(\xi_i)-J_N^*(x) \big |\\ &\leq  \big | J_N^\mathrm{a} \big ( \xi_i ,\D(k_i) \big ) -J_N^*(\xi_i) \big |+2M \norm{x-\xi_i} \to 0.
	\end{align*}
	Thus, we showed that $J_N^\mathrm{a}(x,\D(k))\to J_N^*(x)$ which implies due to \eqref{eq:sandwich} the desired result \eqref{eq:convergenceSW}.
\end{pf}
\begin{figure*}
	\centering
	\resizebox{0.34\linewidth}{!}{
		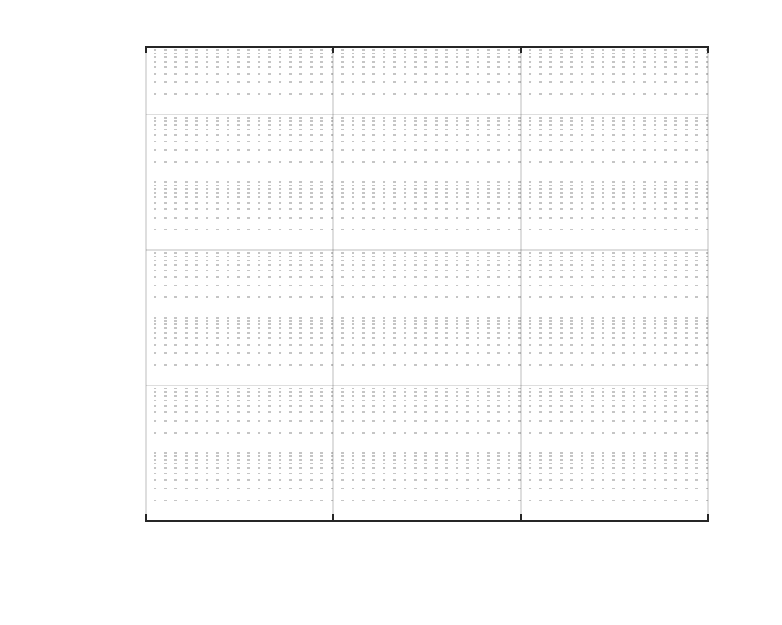}\hspace{-0.5cm}
	\resizebox{0.34\linewidth}{!}{
		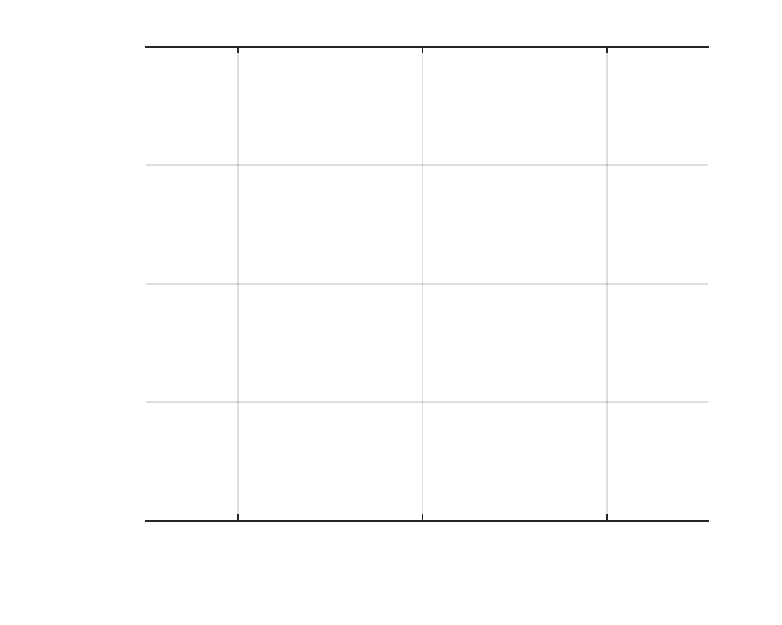}\hspace{-0.5cm}
	\resizebox{0.34\linewidth}{!}{
		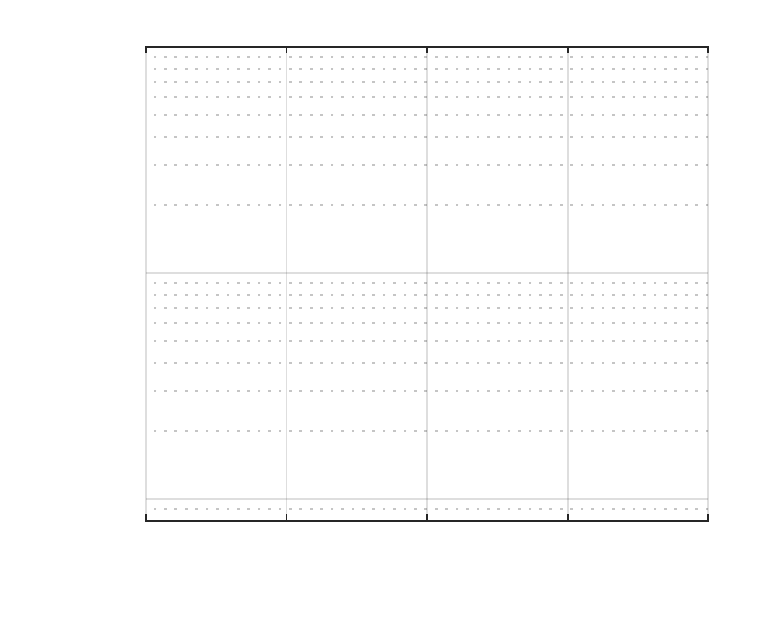}\\[0.3cm]
	\caption{Simulation results of the double integrator example from Section \ref{sec:linAcEx}. \emph{Left:} The suboptimality of the costs for the spatial and temporal warm start $J_N(\Psi_{\mathrm {s/tw}}(k),f(x(k),\Pi_0U(k))-J_N^*(f(x(k),\Pi_0U(k)))$ over time when applying $i_\mathrm{T} (k)=2$ gradient descent steps per optimization. \emph{Middle:} Closed-loop trajectory ({\color{orange}  \textsf{\textbf{\scriptsize o}}}) compared to the original non-learning method~(\raisebox{-0.06cm}{\color{blue} \Large $\cdot$}) starting from the same initial condition. State constraints (dashed) and optimal limit cycle (solid) are also depicted. \emph{Right:} Performance vs. computation time comparison for $i_\mathrm{T}(k)=1,\dots,10$. Averaged over 10 runs ($N_\mathrm{sim}=3000$) with random initial conditions.}\label{fig:acExLin}
\end{figure*}
\begin{rem}\label{rem:initCH}
	The assumption $x\in \mathrm{int}\, \conv \Dx(k_0)$ for some $k_0\in\N$ is rather technical. 
	Still, it is possible that a recurrent point never lies in the interior of $\conv \Dx(k)$ as $k\to\infty$. 
	A straightforward solution to this problem is to initialize the data set $\D(0)$ with some points $(x_i,U_i,J_N(U_i,x_i))$, $i=1,\dots,K$. 
	This solves the issue for all $x\in \interior \conv \{x_1,\dots,x_K\}$ even for arbitrary $U_i$. 
	Hence, a good choice for the points $x_i$ are the extreme points of the polytopic feasible set.
	By adding the point $(0,0,0)$ to $\D(0)$ the prior knowledge about the optimal input at the origin can be included.
\end{rem}
\subsection{Algorithm and implementation}\label{sec:algo}
For online learning, it is crucial to solve \eqref{eq:psi_sw_lin} in an efficient way, since a warm start solution must be provided within one sampling period. 
The convex hull computation gets increasingly burdensome with more data points and higher state space dimensions.
Therefore, we rely on an incremental convex hull algorithm and split the computation in two parts: the convex hull update and the spatial warm start generation. Detailed implementations of both parts can be found in the Appendix \ref{A:algorithm}. 
The update can be performed on a different time scale and does not need to be executed in real time, as is demonstrated in Section \ref{sec:App}.
Notice that the number of data points that is needed to get a good approximation of the MPC control law scales with the dimension of the $\omega$-limit set $\Omega$, which can be significantly lower than the dimension of the state space.
Moreover, the computational complexity of the learning algorithm is independent of the length of the prediction horizon.
Finally, we emphasize that stability is always guaranteed.
Hence, storage and processing restrictions only limit the performance improvement but do not jeopardize the stability of the controller.
\subsection{Example: Double Integrator}\label{sec:linAcEx} 
We demonstrate the learning scheme presented in this section for a numerical example. 
We consider the discrete-time double integrator system from \cite{FELLER2017} of the form \begin{align}
	x(k+1)=\left[\begin{smallmatrix}
	1\ \ & T_\mathrm{s} \\ 0\ \  & 1
	\end{smallmatrix}\right] x(k)+ \left[\begin{smallmatrix}
	T_\mathrm{s}^2 \\T_\mathrm{s}
	\end{smallmatrix}\right]u(k)+w(k)
\end{align}
with sampling time $T_\mathrm{s}=\unit[0.1]{s}$. 
Input and state constraints are inherited from \cite{FELLER2017} $\mathcal U = \{u \in\R : |u| \leq  1\}$ and $\mathcal X = \{x \in\R^2 : -2 \leq  x_1 \leq 3, |x_2| \leq 1\}$, and likewise the parameters of the cost function and of the backtracking line search optimization with gradient descent search direction.
The external disturbance is quasiperiodic $w(k)= 0.09\,[\begin{matrix} \sin(kT_\mathrm{s}) &\cos(kT_\mathrm{s})\end{matrix}]^\top$.

For two optimization iterations $i_\mathrm{T}(k)\equiv i_\mathrm{T}=2$ and a randomly chosen initial condition at $x(0)\approx[0.9\ -0.9]^\top$, we have simulated the closed-loop behavior of the presented MPC scheme for $N_\mathrm{sim}=3000$ time steps. 
In the left subplot of Fig. \ref{fig:acExLin}, the suboptimality of both warm starts, i.e. the difference of their costs to the value function, is depicted over time. 
As was to be expected, the spatial warm start performs first poorly since too few data is available, but then improves over time while more and more data is collected until it significantly outperforms the temporal warm start. 
Further, we can see that the learning rate is limited by the optimization update operator since two gradient descent iterations per time step result in a slow rate of convergence.

In order to highlight the improvements offered by the proposed learning scheme, we compare it to the original anytime MPC \cite{FELLER2017}, i.e. to Algorithm \ref{algo:anyLearn} without lines 4 and 10. 
To this end, the state trajectories of both methods are depicted in the middle subplot of Fig. \ref{fig:acExLin} in the phase portrait until time $k=700$ together with the limit cycle of the optimally controlled system, where the optimization problem \eqref{eq:mpc_prob} is solved with \textsc{Matlab}'s \texttt{fminunc}. 
As we can see, learning significantly improves the performance and leads to much better constraint satisfaction. 
First, both trajectories are identical since the spatial warm start is never used, but as more data is collected the trajectories deviate and the one from the proposed scheme approaches the optimal one.
In particular, the trajectories deviate as soon as they enter the interior of $\conv \Dx(k)$, where for the proposed scheme the convergence result of Theorem \ref{thm:convergence} holds.

Thus far, we have not considered the fact that generating the spatial warm start consumes computation time that might be better invested in doing more optimization iterations with the temporal warm start. 
To investigate this, we run the simulation for different numbers of optimization iterations $i_\mathrm{T}(k)\equiv i_\mathrm{T}=1,\dots ,10$ and for ten different initial conditions. 
For each $i_\mathrm{T}$, we average the computation times and the costs to obtain a point in the right subplot of of Fig. \ref{fig:acExLin}. 
We can see that the proposed learning scheme is way closer to the optimal performance than the original anytime MPC without learning.

%% file: plots/WarmStarts.pdf_tex
\begingroup%
  \makeatletter%
  \providecommand\color[2][]{%
    \errmessage{(Inkscape) Color is used for the text in Inkscape, but the package 'color.sty' is not loaded}%
    \renewcommand\color[2][]{}%
  }%
  \providecommand\transparent[1]{%
    \errmessage{(Inkscape) Transparency is used (non-zero) for the text in Inkscape, but the package 'transparent.sty' is not loaded}%
    \renewcommand\transparent[1]{}%
  }%
  \providecommand\rotatebox[2]{#2}%
  \newcommand*\fsize{\dimexpr\f@size pt\relax}%
  \newcommand*\lineheight[1]{\fontsize{\fsize}{#1\fsize}\selectfont}%
  \ifx\svgwidth\undefined%
    \setlength{\unitlength}{225bp}%
    \ifx\svgscale\undefined%
      \relax%
    \else%
      \setlength{\unitlength}{\unitlength * \real{\svgscale}}%
    \fi%
  \else%
    \setlength{\unitlength}{\svgwidth}%
  \fi%
  \global\let\svgwidth\undefined%
  \global\let\svgscale\undefined%
  \makeatother%
  \begin{picture}(1,0.65)(0, 0.05)%
    \lineheight{1}%
    \setlength\tabcolsep{0pt}%
    \put(0,0){\includegraphics[width=\unitlength,page=1]{plots/WarmStarts.pdf}}%
    \put(0.18666667,0.08355567){\makebox(0,0)[t]{\lineheight{1.25}\smash{\begin{tabular}[t]{c}0\end{tabular}}}}%
    \put(0.42666667,0.08355567){\makebox(0,0)[t]{\lineheight{1.25}\smash{\begin{tabular}[t]{c}100\end{tabular}}}}%
    \put(0.66666667,0.08355567){\makebox(0,0)[t]{\lineheight{1.25}\smash{\begin{tabular}[t]{c}200\end{tabular}}}}%
    \put(0.90666667,0.08355567){\makebox(0,0)[t]{\lineheight{1.25}\smash{\begin{tabular}[t]{c}300\end{tabular}}}}%
    \put(0.546667,0.02422233){\makebox(0,0)[t]{\lineheight{1.25}\smash{\begin{tabular}[t]{c}time $kT_{\mathrm s}$\end{tabular}}}}%
    \put(0,0){\includegraphics[width=\unitlength,page=2]{plots/WarmStarts.pdf}}%
    \put(0.17155567,0.12){\makebox(0,0)[rt]{\lineheight{1.25}\smash{\begin{tabular}[t]{r}$10^{-2}$\end{tabular}}}}%
    \put(0.17155567,0.29333333){\makebox(0,0)[rt]{\lineheight{1.25}\smash{\begin{tabular}[t]{r}$10^{0}$\end{tabular}}}}%
    \put(0.17155567,0.46666667){\makebox(0,0)[rt]{\lineheight{1.25}\smash{\begin{tabular}[t]{r}$10^{2}$\end{tabular}}}}%
    \put(0.17155567,0.64){\makebox(0,0)[rt]{\lineheight{1.25}\smash{\begin{tabular}[t]{r}$10^{4}$\end{tabular}}}}%
    \put(0.07,0.436667){\rotatebox{90}{\makebox(0,0)[t]{\lineheight{1.25}\smash{\begin{tabular}[t]{c} suboptimality $J_N(\cdot)-J_N^*$ \end{tabular}}}}}%
    \put(0,0){\includegraphics[width=\unitlength,page=3]{plots/WarmStarts.pdf}}%
    \put(0.46666667,0.68587433){\makebox(0,0)[lt]{\lineheight{1.25}\smash{\begin{tabular}[t]{l}spatial warm start\end{tabular}}}}%
    \put(0,0){\includegraphics[width=\unitlength,page=4]{plots/WarmStarts.pdf}}%
    \put(0.46666667,0.633459){\makebox(0,0)[lt]{\lineheight{1.25}\smash{\begin{tabular}[t]{l}temporal warm start\end{tabular}}}}%
    \put(0,0){\includegraphics[width=\unitlength,page=5]{plots/WarmStarts.pdf}}%
  \end{picture}%
\endgroup%

%% file: plots/LinAcPhasePortrait.pdf_tex
\begingroup%
  \makeatletter%
  \providecommand\color[2][]{%
    \errmessage{(Inkscape) Color is used for the text in Inkscape, but the package 'color.sty' is not loaded}%
    \renewcommand\color[2][]{}%
  }%
  \providecommand\transparent[1]{%
    \errmessage{(Inkscape) Transparency is used (non-zero) for the text in Inkscape, but the package 'transparent.sty' is not loaded}%
    \renewcommand\transparent[1]{}%
  }%
  \providecommand\rotatebox[2]{#2}%
  \newcommand*\fsize{\dimexpr\f@size pt\relax}%
  \newcommand*\lineheight[1]{\fontsize{\fsize}{#1\fsize}\selectfont}%
  \ifx\svgwidth\undefined%
    \setlength{\unitlength}{225bp}%
    \ifx\svgscale\undefined%
      \relax%
    \else%
      \setlength{\unitlength}{\unitlength * \real{\svgscale}}%
    \fi%
  \else%
    \setlength{\unitlength}{\svgwidth}%
  \fi%
  \global\let\svgwidth\undefined%
  \global\let\svgscale\undefined%
  \makeatother%
  \begin{picture}(1,0.65)(0, 0.05)%
    \lineheight{1}%
    \setlength\tabcolsep{0pt}%
    \put(0,0){\includegraphics[width=\unitlength,page=1]{plots/LinAcPhasePortrait.pdf}}%
    \put(0.30469933,0.08355567){\makebox(0,0)[t]{\lineheight{1.25}\smash{\begin{tabular}[t]{c}-2\end{tabular}}}}%
    \put(0.540765,0.08355567){\makebox(0,0)[t]{\lineheight{1.25}\smash{\begin{tabular}[t]{c}0\end{tabular}}}}%
    \put(0.77683067,0.08355567){\makebox(0,0)[t]{\lineheight{1.25}\smash{\begin{tabular}[t]{c}2\end{tabular}}}}%
    \put(0.546667,0.02266667){\makebox(0,0)[t]{\lineheight{1.25}\smash{\begin{tabular}[t]{c}state $x_1$\end{tabular}}}}%
    \put(0,0){\includegraphics[width=\unitlength,page=2]{plots/LinAcPhasePortrait.pdf}}%
    \put(0.17155567,0.12){\makebox(0,0)[rt]{\lineheight{1.25}\smash{\begin{tabular}[t]{r}-2\end{tabular}}}}%
    \put(0.17155567,0.27166667){\makebox(0,0)[rt]{\lineheight{1.25}\smash{\begin{tabular}[t]{r}-1\end{tabular}}}}%
    \put(0.17155567,0.42333333){\makebox(0,0)[rt]{\lineheight{1.25}\smash{\begin{tabular}[t]{r}0\end{tabular}}}}%
    \put(0.17155567,0.575){\makebox(0,0)[rt]{\lineheight{1.25}\smash{\begin{tabular}[t]{r}1\end{tabular}}}}%
    \put(0.17155567,0.72666667){\makebox(0,0)[rt]{\lineheight{1.25}\smash{\begin{tabular}[t]{r}2\end{tabular}}}}%
    \put(0.1,0.436667){\rotatebox{90}{\makebox(0,0)[t]{\lineheight{1.25}\smash{\begin{tabular}[t]{c}state $x_2$\end{tabular}}}}}%
    \put(0,0){\includegraphics[width=\unitlength,page=3]{plots/LinAcPhasePortrait.pdf}}%
  \end{picture}%
\endgroup%

%% file: plots/timeVperf.pdf_tex
\begingroup%
  \makeatletter%
  \providecommand\color[2][]{%
    \errmessage{(Inkscape) Color is used for the text in Inkscape, but the package 'color.sty' is not loaded}%
    \renewcommand\color[2][]{}%
  }%
  \providecommand\transparent[1]{%
    \errmessage{(Inkscape) Transparency is used (non-zero) for the text in Inkscape, but the package 'transparent.sty' is not loaded}%
    \renewcommand\transparent[1]{}%
  }%
  \providecommand\rotatebox[2]{#2}%
  \newcommand*\fsize{\dimexpr\f@size pt\relax}%
  \newcommand*\lineheight[1]{\fontsize{\fsize}{#1\fsize}\selectfont}%
  \ifx\svgwidth\undefined%
    \setlength{\unitlength}{225bp}%
    \ifx\svgscale\undefined%
      \relax%
    \else%
      \setlength{\unitlength}{\unitlength * \real{\svgscale}}%
    \fi%
  \else%
    \setlength{\unitlength}{\svgwidth}%
  \fi%
  \global\let\svgwidth\undefined%
  \global\let\svgscale\undefined%
  \makeatother%
  \begin{picture}(1,0.65)(0, 0.05)%
    \lineheight{1}%
    \setlength{\fboxsep}{0.5pt}%
    \setlength\tabcolsep{0pt}%
    \put(0,0){\includegraphics[width=\unitlength,page=1]{plots/timeVperf.pdf}}%
    \put(0.18666667,0.08355567){\makebox(0,0)[t]{\lineheight{1.25}\smash{\begin{tabular}[t]{c}0\end{tabular}}}}%
    \put(0.36666667,0.08355567){\makebox(0,0)[t]{\lineheight{1.25}\smash{\begin{tabular}[t]{c}0.01\end{tabular}}}}%
    \put(0.54666667,0.08355567){\makebox(0,0)[t]{\lineheight{1.25}\smash{\begin{tabular}[t]{c}0.02\end{tabular}}}}%
    \put(0.72666667,0.08355567){\makebox(0,0)[t]{\lineheight{1.25}\smash{\begin{tabular}[t]{c}0.03\end{tabular}}}}%
    \put(0.90666667,0.08355567){\makebox(0,0)[t]{\lineheight{1.25}\smash{\begin{tabular}[t]{c}0.04\end{tabular}}}}%
    \put(0.546667,0.02422233){\makebox(0,0)[t]{\lineheight{1.25}\smash{\begin{tabular}[t]{c} average computation time per time step\end{tabular}}}}%
    \put(0,0){\includegraphics[width=\unitlength,page=2]{plots/timeVperf.pdf}}%
    \put(0.17155567,0.14803733){\makebox(0,0)[rt]{\lineheight{1.25}\smash{\begin{tabular}[t]{r}$10^{2}$\end{tabular}}}}%
    \put(0.17155567,0.437352){\makebox(0,0)[rt]{\lineheight{1.25}\smash{\begin{tabular}[t]{r}$10^{3}$\end{tabular}}}}%
    \put(0.17155567,0.72666667){\makebox(0,0)[rt]{\lineheight{1.25}\smash{\begin{tabular}[t]{r}$10^{4}$\end{tabular}}}}%
    \put(0.075,0.436667){\rotatebox{90}{\makebox(0,0)[t]{\lineheight{1.25}\smash{\begin{tabular}[t]{c}average cost per run\end{tabular}}}}}%
    \put(0,0){\includegraphics[width=\unitlength,page=3]{plots/timeVperf.pdf}}%
    \put(0.52333333,0.66633333){\makebox(0,0)[lt]{\lineheight{1.25}\smash{\begin{tabular}[t]{l}\ without learning \end{tabular}}}}%
    \put(0,0){\includegraphics[width=\unitlength,page=4]{plots/timeVperf.pdf}}%
    \put(0.52333333,0.61466667){\makebox(0,0)[lt]{\lineheight{1.25}\smash{\begin{tabular}[t]{l}\ with learning\end{tabular}}}}%
    \put(0,0){\includegraphics[width=\unitlength,page=5]{plots/timeVperf.pdf}}%
    \put(0.52333333,0.563){\makebox(0,0)[lt]{\lineheight{1.25}\smash{\begin{tabular}[t]{l}\ optimal\end{tabular}}}}%
    \put(0,0){\includegraphics[width=\unitlength,page=6]{plots/timeVperf.pdf}}%
    \put(0.23,0.32){\small \fcolorbox{white}{white}{$i_\mathrm{T}=1$}}%
    \put(0.2,0.67){\small \fcolorbox{white}{white}{$i_\mathrm{T}=1$}}%
    \put(0.75,0.32){\small \fcolorbox{white}{white}{$i_\mathrm{T}=10$}}%
    \put(0.79,0.2){\small  \fcolorbox{white}{white}{$i_\mathrm{T}=10$}}%
  \end{picture}%
\endgroup%

%% file: 4_nonlin.tex
\section{Leveraging data in real-time nonlinear MPC}\label{sec:nonLinMPC}
In the general case of the nonlinear MPC scheme \eqref{eq:sysx}--\eqref{eq:D}, \eqref{eq:psi_w}, we cannot expect that the cost function $J_N$ is convex. Therefore, we propose in this section a different learning method that does not depend on convexity, but instead exploits Lipschitz continuity.
\subsection{A spatial warm start based on Lipschitz continuity}
The idea for the spatial warm start generation for $\bar x\in\R^n$ is to find an input sequence from the data $(U,x,J)\in\D$ such that $x$ is close to $\bar x$ and $J$ is small. This directly leads to the following spatial warm start rule
\begin{subequations}\label{eq:psi_sw_nonl}
	\begin{align}\nonumber
	\big (\Psi_\mathrm{sw}(\bar x,\mathcal D),x^*&(\bar x,\D),J^*(\bar x,\D) \big )\\ \label{eq:xJstar} &=\argmin_{(U,x,J)\in \mathcal D} J+L(U)\norm{x-\bar x}\\
	J_N^\mathrm{a}(\bar x,\mathcal D) &= \min_{(U,x,J)\in \mathcal D} J+L(U)\norm{x-\bar x}\label{eq:JNagenau}
	\end{align}
\end{subequations}
where $L(U)\in\Rpp$ is a Lipschitz parameter that might depend on $U$, and $J_N^\mathrm{a} (\cdot ,\D)$ is an approximation of the value function $J_N^*$ based on the data $\D$. An illustration of the spatial warm start \eqref{eq:psi_sw_nonl} is given in Fig. \ref{fig:psi_sw_nonl}. If we assume $J_N(U,\cdot)$ to be Lipschitz continuous with constant $L(U)$, then $J_N^\mathrm{a}(\cdot,\D)$ upper bounds the value function and the spatial warm start as shown in the following Lemma. Similar ideas to approximate an unknown function based on Lipschitz continuity from sampling data have been proposed in \cite{BELIAKOV2006}, \cite{CALLIESS2014}, and \cite{ZABINSKY2003}.

\begin{assum}\label{assum:lipschitz}
	Assume that $J_N(U, \cdot)$ is Lipschitz continuous with constant $L(U)$, i.e. for all $x,y \in\R^n$,
	\begin{align}\label{eq:Lipsch}
	\big| J_N(U,x)-J_N(U,y)\big |\leq L(U) \norm{x-y}.
	\end{align}
\end{assum}
\begin{lem}\label{lem:JNaPropsNL}
	Let Assumption \ref{assum:lipschitz} hold, then 
	\begin{align}\label{eq:JNaUpb}
	J_N\big (\Psi_\mathrm{sw}(x,\mathcal D),x\big )\leq J_N^\mathrm{a} (x,\mathcal D)
	\end{align}
	and $J_N^\mathrm{a}(\cdot,\D)$ is also Lipschitz continuous with Lipschitz constant $L_\D=\max_{U\in\Du} L(U)$.
\end{lem}
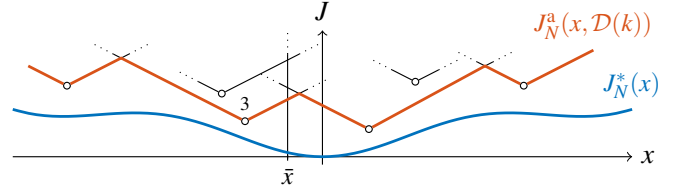
\begin{figure}
	\centering
	\resizebox{\linewidth}{!}{
		\begin{tikzpicture}[xscale=0.98,yscale=0.5]
		\draw[->] (-4,0)--(4,0) node [right] {$x$};
		\draw[->] (0,-0.2)--(0,3.2) node [above] {$J$};
		\node [inner sep=0.03cm, circle, draw] (e3) at (-3.3,1.8) {};
		\node [inner sep=0.03cm, circle, draw] (e6) at (-1,0.9) {};
		\node [inner sep=0.03cm, circle, draw] (e2) at (0.6,0.7) {};
		\node [inner sep=0.03cm, circle, draw] (e1) at (2.6,1.8) {};
		\node [inner sep=0.03cm, circle, draw] (e5) at (1.2,1.9) {};
		\node [inner sep=0.03cm, circle, draw] (e4) at (-1.3,1.6) {};
		\draw (e4) -- (-0.3,2.6);
		\draw [dotted] (-0.3,2.6) -- (-0.05,2.85);
		\draw (e4) -- (-1.6,1.9);
		\draw [dotted] (-1.6,1.9) -- (-1.85,2.15);
		\draw (-0.45,-0.1) -- (-0.45,2.6);
		\draw [dotted]  (-0.45,2.6) -- (-0.45,3.15);
		\draw (-0.55,1.85) -- (-0.3,1.6);
		\draw [dotted] (-0.55,1.85) -- (-0.8,2.1);
		\draw (-2.6,2.5) -- (-2.5,2.6);
		\draw [dotted] (-2.5,2.6) -- (-2.25,2.85);
		\draw (-2.6,2.5) -- (-2.7,2.6);
		\draw [dotted] (-2.7,2.6) -- (-2.95,2.85);
		\draw (e5) -- (0.9,2.2);
		\draw [dotted] (0.9,2.2) -- (0.65,2.45);
		\draw (e5) -- +(0.3,0.3);
		\draw [dotted] (e5)+(0.3,0.3) -- +(0.55,0.55);
		\draw (2.1,2.3) -- (2,2.4);
		\draw [dotted] (2,2.4) -- (1.75,2.65);
		\draw (2.1,2.3) -- (2.2,2.4);
		\draw [dotted] (2.2,2.4) -- (2.45,2.65);
		\draw (-0.3,1.6) -- (-0.2,1.7);
		\draw [dotted] (-0.2,1.7) -- (0.05,1.95);
		\draw (-0.45,-0.05)node [below] {\small $\bar x$} -- (-0.45,0.05);
		\draw [smooth, domain=-4:4, samples=40,variable=\x,blue,line width=1pt] plot ({\x},{5/6*(0.09*(\x*\x)+0.9*sin(abs(\x)*45)*sin(abs(\x)*45))}) node [above]{\small $J_N^*(x)$};
		\draw[orange,line width =1pt] (-3.8,2.3) -- (e3) -- (-2.6,2.5) -- (e6)node [above] {\color{black}\scriptsize $3$} -- (-0.3,1.6) -- (e2) -- (2.1,2.3) -- (e1) -- (3.5,2.7) node [above] {\small $J^{\mathrm a}_N (x,\mathcal D(k))$};
		\end{tikzpicture}}
	\caption{Illustration of the spatial warm start \eqref{eq:psi_sw_nonl} in the $(x,J)$-domain with data points $(x_i,J_i)\in \Dxj(k)$ ({\scriptsize$\circ$}). 
		For each point, the cone $J_i+L(U_i)\norm{x-x_i}$ is indicated as well as the minimum over all $i$ for each $x$ given by $J_N^\mathrm{a}(x,\D(k))$, which is an upper bound for the value function. 
		The spatial warm start solution at $\bar x$ is the input from point 3.
}\label{fig:psi_sw_nonl}
\end{figure}
\begin{figure*}[t!]
\centering
\resizebox{0.5\linewidth}{!}{
	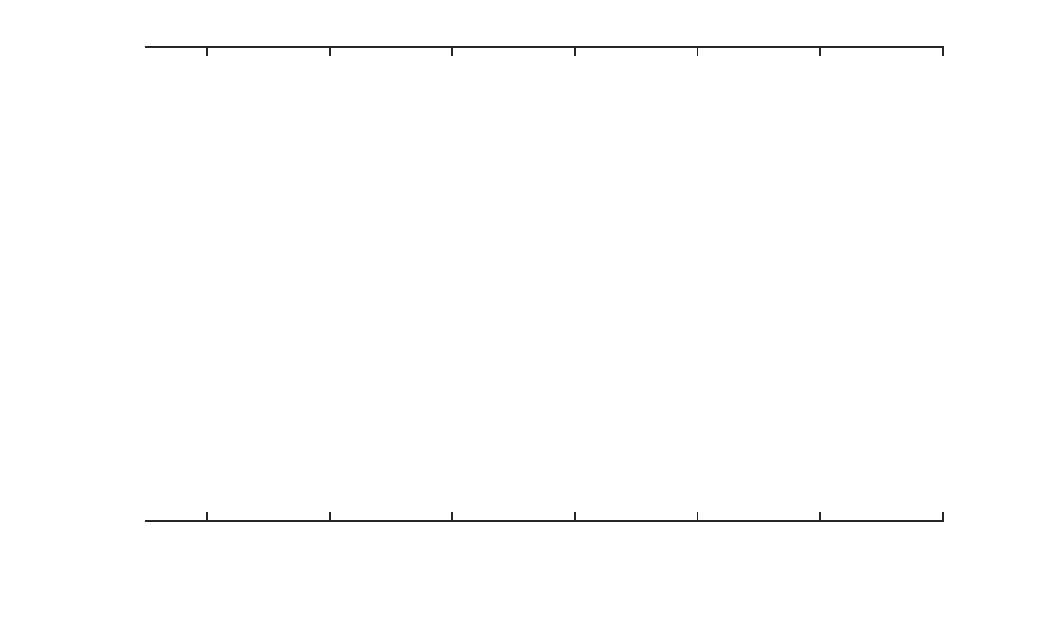}\hspace{-0.1cm}
\resizebox{0.5\linewidth}{!}{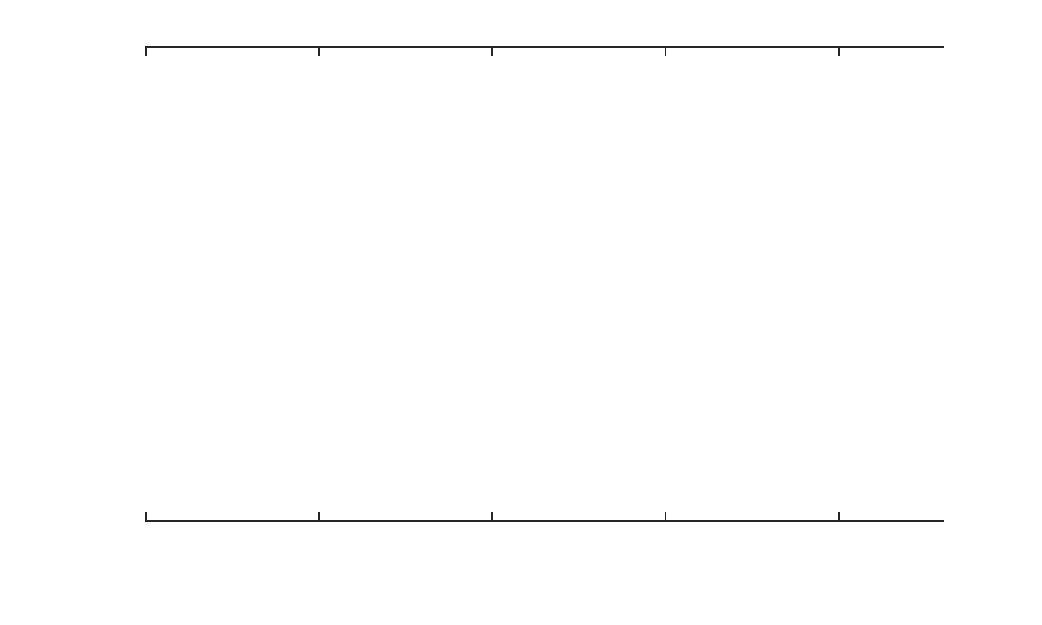}\\[0.1cm]
\caption{Simulation results of the unicycle example from Section \ref{sec:nonlAcEx} controlled by the proposed MPC scheme for nonlinear systems. \emph{Left:} State trajectory for several runs and the optimal trajectory (computed with \textsc{Matlab}'s \texttt{patternsearch}). \emph{Right:} Accumulated cost per run over the number of the runs. Additionally, the optimal cost is depicted.}\label{fig:acExNonl}
\end{figure*}
For a proof, see Appendix \ref{sec:appendix_proofs}.
With Lemma \ref{lem:JNaPropsNL}, we obtain a convergence result similar to Theorem \ref{thm:convergence}.
\begin{thm}\label{thm:convergence_nonl}
	Consider the nonlinear MPC scheme presented in \eqref{eq:sysx}--\eqref{eq:omega}, \eqref{eq:psi_w} with the spatial warm start operator \eqref{eq:psi_sw_nonl}, further let Assumption \ref{assum:lipschitz} hold, and let $\forall k\geq 1:L(U(k))\leq M\in\R$ be upper bounded. Then for all $x\in \Omega$
	\begin{align}
	\lim_{k\to \infty} J_N \big (\Psi_\mathrm{sw}\big (x,\mathcal D(k) \big ),x\big )=J_N^*(x).
	 \end{align}
\end{thm}
\begin{pf} Step 1) $J_N^\mathrm{a}(x,\D(k))$ converges:
In view of Lemma \ref{lem:JNaPropsNL}, we have $
	J_N^*(x)\leq J_N(\Psi_{\mathrm{sw}}(x,\mathcal D(k)),x)\leq J_N^{\mathrm{a}}(x,\mathcal D(k)).$ 
	Thus by showing $J_N^\mathrm{a} (x,\D(k))\to J_N^*(x)$ for $x\in\Omega$ we will prove the theorem. We also see from this inequality that $J_N^\mathrm{a}(x,\D(k))$ is bounded from below by $J_N^*(x)$. Further $J_N^\mathrm{a}(x,\D(k))$ is decreasing since the minimum over a larger set $\D(k+1)\supseteq \D(k)$ is smaller or equal and thus $J_N^\mathrm{a}(x,\D(k))$ must converge to some value. \\
	Step 2) $J_N^\mathrm{a}(x,\D(k))$ converges to $J_N^*(x)$: This step is analogous to step 2) in the proof of Theorem \ref{thm:convergence} and thus moved to the Appendix \ref{sec:appendix_proofs}.
\end{pf}
\begin{rem}
	Even tough Assumption \ref{assum:lipschitz} might be restrictive, if we assume that $x(k)$ does not grow unbounded, but stays in some compact region $\mathcal C_\mathrm{x}\subset \R^{n}$, then the Lipschitz constants $L(U(k))$ do not need to apply globally, but only on $\mathcal C_\mathrm{x}$ in order to obtain the same result. 
	If further $U(k)$ stays in some compact set $\mathcal C_\mathrm{U}\subseteq\R^{Nm}$, then the assumption $L(U(k)) \leq M$ is also satisfied with $M=\max_{U\in \mathcal C_\mathrm{U}} L(U)<\infty$ if $J_N$ is continuous in $U$. 
\end{rem}
\begin{rem}
	It can be quite challenging to satisfy \eqref{eq:classKfncts} and \eqref{eq:ineq} for general nonlinear MPC algorithms. 
	As often done in practice, one can still implement the learning scheme with spatial warm start \eqref{eq:psi_sw_nonl}. 
	The controller performance will improve as long as the optimization iteration does reduce the costs, even if \eqref{eq:classKfncts} and \eqref{eq:ineq} do not hold. 
	However, a prior stability guarantee depends on the nonlinear iteration scheme.
\end{rem}
\begin{rem}
	If in a specific setup the cost function $J_N$ is nevertheless known to be convex, then the learning scheme of Section \ref{sec:linMPC} can be used and probably leads to a better upper bound $J_N^\mathrm{a}$; for example, see Fig. \ref{fig:psi_sw_nonl} where the convex hull over the data points would result in a lower upper bound $J_N^\mathrm{a}(\cdot ,\D)$ inside $\conv \D_x$. 
	On the other hand, if one faces a linear MPC problem with a nonconvex cost function $J_N$, due to nonconvex stage or terminal costs, then the learning scheme from this section can be applied.
\end{rem}

\subsection{Example: Unicycle}\label{sec:nonlAcEx}
In this section, we implement the nonlinear real-time MPC learning scheme for the unicycle model
\begin{align}
\dot x = f_\mathrm{c}(x,u)= \begin{bmatrix}
u_1 \cos(x_3) &u_1 \sin(x_3)& u_2
\end{bmatrix}^\top
\end{align}
with state and input vectors $x=[x_1\ x_2\ x_3]$, $u=[ u_1 \  u_2]$ respectively, as well as the sampling time $T_\mathrm{s} =\unit[0.1]{s}$. 
The control task is to drive the unicycle to $(x_1,x_2)=(0,0)$ facing into positive $x_1$ direction, i.e. $x_3=2z\pi$ for some $z\in\Z$ while keeping $\norm{u}_\infty<1$. 
Therefore, the stage cost in \eqref{eq:mpc_prob} is chosen as $l(x,u)=0.1 \sin\left (\frac{x_3} 2 \right )^2+\sqrt[4]{1+x_1^2+x_2^2}-1+u_1^8+u_2^8$ to ensure positive definiteness with respect to $u=0$, $(x_1,x_2,x_3)=(0,0,2z\pi)$, prevent $x_1^2+x_2^2$ from getting too steep, while ensuring differentiability at the origin, and let $u$ shoot up as $\norm{u}_\infty>1$.
The terminal cost is set to $F(x)=100\,l(x,0)$.
For this example, the learning scheme from Section \ref{sec:linMPC} cannot be applied since the resulting cost function $J_N$ is not convex.
However, the costs turn out to be Lipschitz continuous in $x$ with Lipschitz constant $L(U)$ computed in Appendix \ref{A:lipsch}.

We choose $\Psi_{\mathrm {tw}}(U,x)=[U_1^\top,\dots, U_{N-1}^\top, 0^\top]^\top$ as the temporal warm start and the same optimization operator \eqref{eq:psi_o} as for the linear example in Section \ref{sec:linAcEx}, which consists of the gradient descent search direction with backtracking line search. 
Further, we use a constant number of optimization iterations $i_\mathrm{T}(k)=2$. 
We consider iterative learning in a repetitive and fast process, which is a common task in iterative learning control (see \cite{WANG2009}). 
The system is repeatedly started from the same initial condition $x_0=[1,1,1+\pi/2]^\top$, where each run takes $120$ time steps.
Hence, throughout a run, $w(k)$ is zero and after every $120$ time steps it sets $x(k)$ back to $x_0$.

\begin{figure*}
	\centering
	\resizebox{0.45\linewidth}{!}{
		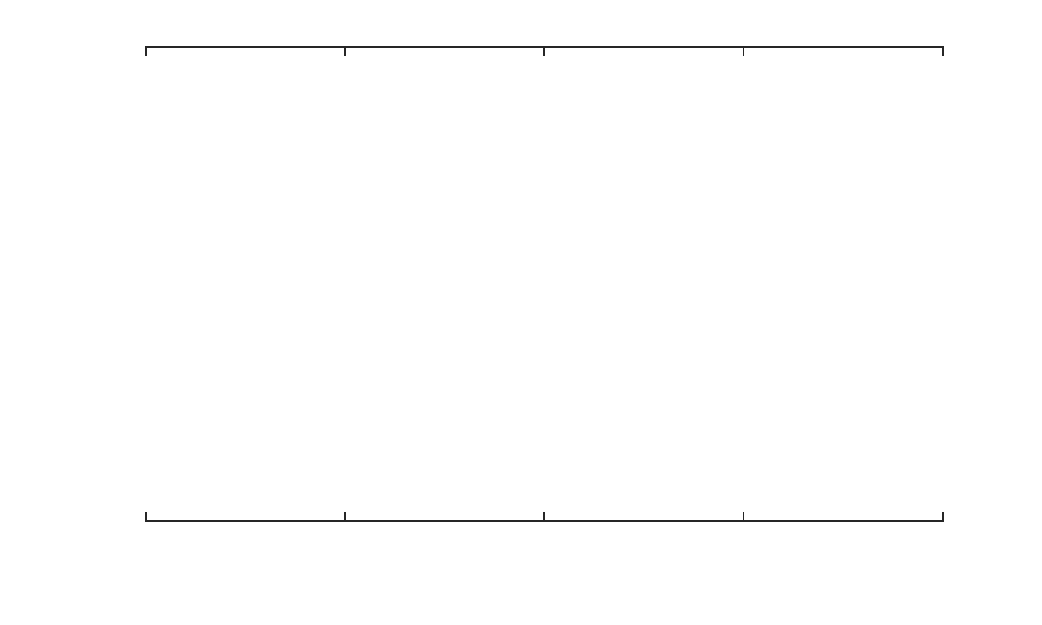}\hspace{-0.5cm}
	\resizebox{0.29\linewidth}{!}{
		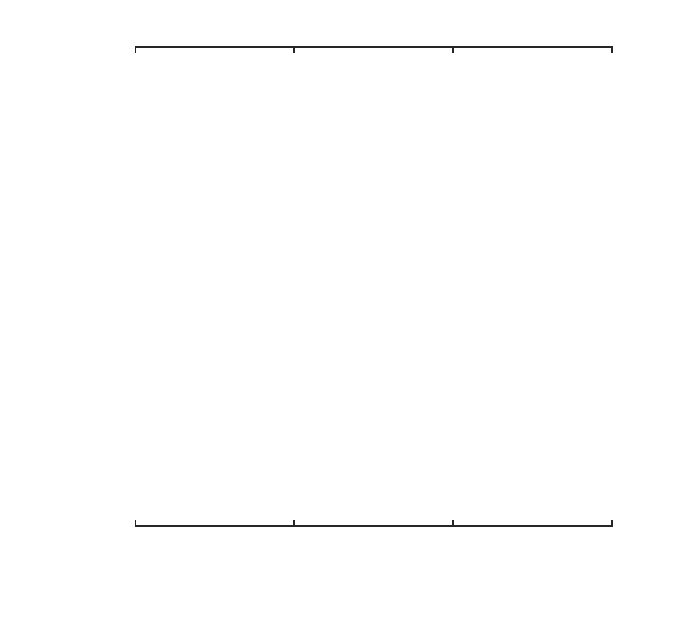}\hspace{-0.5cm}
	\resizebox{0.29\linewidth}{!}{
		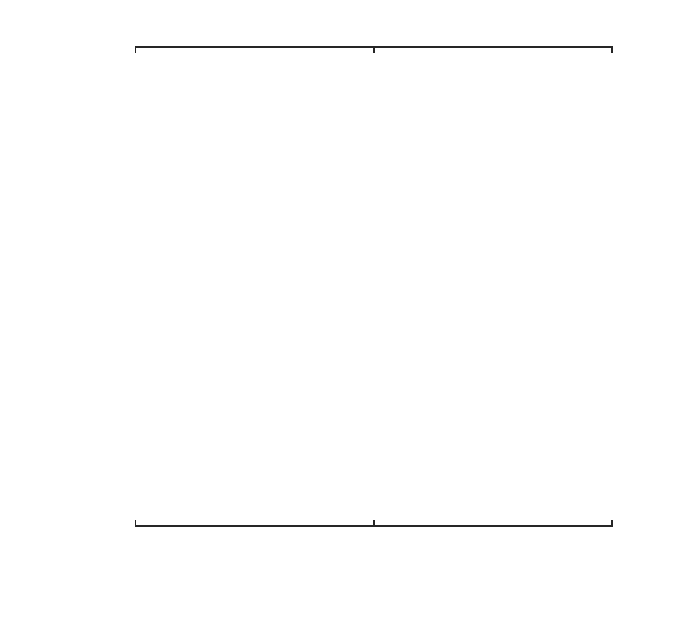}
	\caption{Simulation results of the servomechanism in Section \ref{sec:App}.
		\emph{Left:} The effect of learning is clearly visible in the improved reference tracking behavior over increasing number of periods. 
		\emph{Middle:} The mean tracking error per period of the linear anytime MPC with and without the proposed learning scheme as well as the optimal mean error. 
		\emph{Right:} First ten seconds of the input signal for different periods, as well as the optimal input and the input constraint.}\label{fig:applEx}
\end{figure*}

Although $\Psi_\mathrm{tw}$ and $\Psi_\mathrm{o}$ might not satisfy \eqref{eq:ineq}, the closed loop performance is still satisfying as we can see in the simulation results depicted in Fig. \ref{fig:acExNonl}. 
In the plot, we can see that the optimal behavior is learned within a few runs.
After five runs the unicycle has learned to start driving backwards, and after $20$ runs it has learned to approach the destination driving forwards and is almost indistinguishable from the optimal trajectory.

%% file: plots/NonlPhasePortrait.pdf_tex
\begingroup%
  \makeatletter%
  \providecommand\color[2][]{%
    \errmessage{(Inkscape) Color is used for the text in Inkscape, but the package 'color.sty' is not loaded}%
    \renewcommand\color[2][]{}%
  }%
  \providecommand\transparent[1]{%
    \errmessage{(Inkscape) Transparency is used (non-zero) for the text in Inkscape, but the package 'transparent.sty' is not loaded}%
    \renewcommand\transparent[1]{}%
  }%
  \providecommand\rotatebox[2]{#2}%
  \newcommand*\fsize{\dimexpr\f@size pt\relax}%
  \newcommand*\lineheight[1]{\fontsize{\fsize}{#1\fsize}\selectfont}%
  \ifx\svgwidth\undefined%
    \setlength{\unitlength}{300bp}%
    \ifx\svgscale\undefined%
      \relax%
    \else%
      \setlength{\unitlength}{\unitlength * \real{\svgscale}}%
    \fi%
  \else%
    \setlength{\unitlength}{\svgwidth}%
  \fi%
  \global\let\svgwidth\undefined%
  \global\let\svgscale\undefined%
  \makeatother%
  \begin{picture}(1,0.55)(0, 0.02)%
    \lineheight{1}%
    \setlength\tabcolsep{0pt}%
    \put(0,0){\includegraphics[width=\unitlength,page=1]{plots/NonlPhasePortrait.pdf}}%
    \put(0.19884625,0.06266675){\makebox(0,0)[t]{\lineheight{1.25}\smash{\begin{tabular}[t]{c}0\end{tabular}}}}%
    \put(0.3165385,0.06266675){\makebox(0,0)[t]{\lineheight{1.25}\smash{\begin{tabular}[t]{c}0.2\end{tabular}}}}%
    \put(0.43423075,0.06266675){\makebox(0,0)[t]{\lineheight{1.25}\smash{\begin{tabular}[t]{c}0.4\end{tabular}}}}%
    \put(0.551923,0.06266675){\makebox(0,0)[t]{\lineheight{1.25}\smash{\begin{tabular}[t]{c}0.6\end{tabular}}}}%
    \put(0.6696155,0.06266675){\makebox(0,0)[t]{\lineheight{1.25}\smash{\begin{tabular}[t]{c}0.8\end{tabular}}}}%
    \put(0.78730775,0.06266675){\makebox(0,0)[t]{\lineheight{1.25}\smash{\begin{tabular}[t]{c}1\end{tabular}}}}%
    \put(0.905,0.06266675){\makebox(0,0)[t]{\lineheight{1.25}\smash{\begin{tabular}[t]{c}1.2\end{tabular}}}}%
    \put(0.5225005,0.017){\makebox(0,0)[t]{\lineheight{1.25}\smash{\begin{tabular}[t]{c}state $x_1$\end{tabular}}}}%
    \put(0,0){\includegraphics[width=\unitlength,page=2]{plots/NonlPhasePortrait.pdf}}%
    \put(0.12866675,0.146875){\makebox(0,0)[rt]{\lineheight{1.25}\smash{\begin{tabular}[t]{r}0\end{tabular}}}}%
    \put(0.12866675,0.22270825){\makebox(0,0)[rt]{\lineheight{1.25}\smash{\begin{tabular}[t]{r}0.2\end{tabular}}}}%
    \put(0.12866675,0.29854175){\makebox(0,0)[rt]{\lineheight{1.25}\smash{\begin{tabular}[t]{r}0.4\end{tabular}}}}%
    \put(0.12866675,0.374375){\makebox(0,0)[rt]{\lineheight{1.25}\smash{\begin{tabular}[t]{r}0.6\end{tabular}}}}%
    \put(0.12866675,0.45020825){\makebox(0,0)[rt]{\lineheight{1.25}\smash{\begin{tabular}[t]{r}0.8\end{tabular}}}}%
    \put(0.12866675,0.52604175){\makebox(0,0)[rt]{\lineheight{1.25}\smash{\begin{tabular}[t]{r}1\end{tabular}}}}%
    \put(0.037,0.32750025){\rotatebox{90}{\makebox(0,0)[t]{\lineheight{1.25}\smash{\begin{tabular}[t]{c}state $x_2$\end{tabular}}}}}%
    \put(0,0){\includegraphics[width=\unitlength,page=3]{plots/NonlPhasePortrait.pdf}}%
    \put(0.63,0.329873){\makebox(0,0)[lt]{\lineheight{1.25}\smash{\begin{tabular}[t]{l}\ 1st run\end{tabular}}}}%
    \put(0,0){\includegraphics[width=\unitlength,page=4]{plots/NonlPhasePortrait.pdf}}%
    \put(0.63,0.29132375){\makebox(0,0)[lt]{\lineheight{1.25}\smash{\begin{tabular}[t]{l}\ 5th run\end{tabular}}}}%
    \put(0,0){\includegraphics[width=\unitlength,page=5]{plots/NonlPhasePortrait.pdf}}%
    \put(0.63,0.2527745){\makebox(0,0)[lt]{\lineheight{1.25}\smash{\begin{tabular}[t]{l}\ 10th run\end{tabular}}}}%
    \put(0,0){\includegraphics[width=\unitlength,page=6]{plots/NonlPhasePortrait.pdf}}%
    \put(0.63,0.2142255){\makebox(0,0)[lt]{\lineheight{1.25}\smash{\begin{tabular}[t]{l}\ 20th run\end{tabular}}}}%
    \put(0,0){\includegraphics[width=\unitlength,page=7]{plots/NonlPhasePortrait.pdf}}%
    \put(0.63,0.17567625){\makebox(0,0)[lt]{\lineheight{1.25}\smash{\begin{tabular}[t]{l}\ opt. trajectory\end{tabular}}}}%
    \put(0,0){\includegraphics[width=\unitlength,page=8]{plots/NonlPhasePortrait.pdf}}%
    \put(0.63,0.137127){\makebox(0,0)[lt]{\lineheight{1.25}\smash{\begin{tabular}[t]{l}\ initial point \end{tabular}}}}%
    \put(0,0){\includegraphics[width=\unitlength,page=9]{plots/NonlPhasePortrait.pdf}}%
  \end{picture}%
\endgroup%

%% file: plots/NonlOpt.pdf_tex
\begingroup%
  \makeatletter%
  \providecommand\color[2][]{%
    \errmessage{(Inkscape) Color is used for the text in Inkscape, but the package 'color.sty' is not loaded}%
    \renewcommand\color[2][]{}%
  }%
  \providecommand\transparent[1]{%
    \errmessage{(Inkscape) Transparency is used (non-zero) for the text in Inkscape, but the package 'transparent.sty' is not loaded}%
    \renewcommand\transparent[1]{}%
  }%
  \providecommand\rotatebox[2]{#2}%
  \newcommand*\fsize{\dimexpr\f@size pt\relax}%
  \newcommand*\lineheight[1]{\fontsize{\fsize}{#1\fsize}\selectfont}%
  \ifx\svgwidth\undefined%
    \setlength{\unitlength}{300bp}%
    \ifx\svgscale\undefined%
      \relax%
    \else%
      \setlength{\unitlength}{\unitlength * \real{\svgscale}}%
    \fi%
  \else%
    \setlength{\unitlength}{\svgwidth}%
  \fi%
  \global\let\svgwidth\undefined%
  \global\let\svgscale\undefined%
  \makeatother%
  \begin{picture}(1,0.55)(0, 0.02)%
    \lineheight{1}%
    \setlength\tabcolsep{0pt}%
    \put(0,0){\includegraphics[width=\unitlength,page=1]{plots/NonlOpt.pdf}}%
    \put(0.14,0.06266675){\makebox(0,0)[t]{\lineheight{1.25}\smash{\begin{tabular}[t]{c}0\end{tabular}}}}%
    \put(0.30630425,0.06266675){\makebox(0,0)[t]{\lineheight{1.25}\smash{\begin{tabular}[t]{c}5\end{tabular}}}}%
    \put(0.47260875,0.06266675){\makebox(0,0)[t]{\lineheight{1.25}\smash{\begin{tabular}[t]{c}10\end{tabular}}}}%
    \put(0.638913,0.06266675){\makebox(0,0)[t]{\lineheight{1.25}\smash{\begin{tabular}[t]{c}15\end{tabular}}}}%
    \put(0.80521725,0.06266675){\makebox(0,0)[t]{\lineheight{1.25}\smash{\begin{tabular}[t]{c}20\end{tabular}}}}%
    \put(0.5225005,0.017){\makebox(0,0)[t]{\lineheight{1.25}\smash{\begin{tabular}[t]{c}runs\end{tabular}}}}%
    \put(0,0){\includegraphics[width=\unitlength,page=2]{plots/NonlOpt.pdf}}%
    \put(0.12866675,0.09){\makebox(0,0)[rt]{\lineheight{1.25}\smash{\begin{tabular}[t]{r}250\end{tabular}}}}%
    \put(0.12866675,0.19833325){\makebox(0,0)[rt]{\lineheight{1.25}\smash{\begin{tabular}[t]{r}300\end{tabular}}}}%
    \put(0.12866675,0.30666675){\makebox(0,0)[rt]{\lineheight{1.25}\smash{\begin{tabular}[t]{r}350\end{tabular}}}}%
    \put(0.12866675,0.415){\makebox(0,0)[rt]{\lineheight{1.25}\smash{\begin{tabular}[t]{r}400\end{tabular}}}}%
    \put(0.12866675,0.52333325){\makebox(0,0)[rt]{\lineheight{1.25}\smash{\begin{tabular}[t]{r}450\end{tabular}}}}%
    \put(0.037,0.32750025){\rotatebox{90}{\makebox(0,0)[t]{\lineheight{1.25}\smash{\begin{tabular}[t]{c}cost\end{tabular}}}}}%
    \put(0,0){\includegraphics[width=\unitlength,page=3]{plots/NonlOpt.pdf}}%
    \put(0.645,0.49940575){\makebox(0,0)[lt]{\lineheight{1.25}\smash{\begin{tabular}[t]{l}cost per run\end{tabular}}}}%
    \put(0,0){\includegraphics[width=\unitlength,page=4]{plots/NonlOpt.pdf}}%
    \put(0.645,0.46009425){\makebox(0,0)[lt]{\lineheight{1.25}\smash{\begin{tabular}[t]{l}optimal cost \end{tabular}}}}%
    \put(0,0){\includegraphics[width=\unitlength,page=5]{plots/NonlOpt.pdf}}%
  \end{picture}%
\endgroup%

%% file: plots/ApplExRefTrack.pdf_tex
\begingroup%
  \makeatletter%
  \providecommand\color[2][]{%
    \errmessage{(Inkscape) Color is used for the text in Inkscape, but the package 'color.sty' is not loaded}%
    \renewcommand\color[2][]{}%
  }%
  \providecommand\transparent[1]{%
    \errmessage{(Inkscape) Transparency is used (non-zero) for the text in Inkscape, but the package 'transparent.sty' is not loaded}%
    \renewcommand\transparent[1]{}%
  }%
  \providecommand\rotatebox[2]{#2}%
  \newcommand*\fsize{\dimexpr\f@size pt\relax}%
  \newcommand*\lineheight[1]{\fontsize{\fsize}{#1\fsize}\selectfont}%
  \ifx\svgwidth\undefined%
    \setlength{\unitlength}{300bp}%
    \ifx\svgscale\undefined%
      \relax%
    \else%
      \setlength{\unitlength}{\unitlength * \real{\svgscale}}%
    \fi%
  \else%
    \setlength{\unitlength}{\svgwidth}%
  \fi%
  \global\let\svgwidth\undefined%
  \global\let\svgscale\undefined%
  \makeatother%
  \begin{picture}(1,0.6)%
    \lineheight{1}%
    \setlength\tabcolsep{0pt}%
    \put(0,0){\includegraphics[width=\unitlength,page=1]{plots/ApplExRefTrack.pdf}}%
    \put(0.14,0.06266675){\makebox(0,0)[t]{\lineheight{1.25}\smash{\begin{tabular}[t]{c}0\end{tabular}}}}%
    \put(0.33125,0.06266675){\makebox(0,0)[t]{\lineheight{1.25}\smash{\begin{tabular}[t]{c}5\end{tabular}}}}%
    \put(0.5225,0.06266675){\makebox(0,0)[t]{\lineheight{1.25}\smash{\begin{tabular}[t]{c}10\end{tabular}}}}%
    \put(0.71375,0.06266675){\makebox(0,0)[t]{\lineheight{1.25}\smash{\begin{tabular}[t]{c}15\end{tabular}}}}%
    \put(0.905,0.06266675){\makebox(0,0)[t]{\lineheight{1.25}\smash{\begin{tabular}[t]{c}20\end{tabular}}}}%
    \put(0.5225005,0.017){\makebox(0,0)[t]{\lineheight{1.25}\smash{\begin{tabular}[t]{c}time within period in s\end{tabular}}}}%
    \put(0,0){\includegraphics[width=\unitlength,page=2]{plots/ApplExRefTrack.pdf}}%
    \put(0.12866675,0.092113){\makebox(0,0)[rt]{\lineheight{1.25}\smash{\begin{tabular}[t]{r}0\end{tabular}}}}%
    \put(0.12866675,0.179597){\makebox(0,0)[rt]{\lineheight{1.25}\smash{\begin{tabular}[t]{r}0.2\end{tabular}}}}%
    \put(0.12866675,0.26708075){\makebox(0,0)[rt]{\lineheight{1.25}\smash{\begin{tabular}[t]{r}0.4\end{tabular}}}}%
    \put(0.12866675,0.35456475){\makebox(0,0)[rt]{\lineheight{1.25}\smash{\begin{tabular}[t]{r}0.6\end{tabular}}}}%
    \put(0.12866675,0.44204875){\makebox(0,0)[rt]{\lineheight{1.25}\smash{\begin{tabular}[t]{r}0.8\end{tabular}}}}%
    \put(0.12866675,0.5295325){\makebox(0,0)[rt]{\lineheight{1.25}\smash{\begin{tabular}[t]{r}1\end{tabular}}}}%
    \put(0.037,0.32750025){\rotatebox{90}{\makebox(0,0)[t]{\lineheight{1.25}\smash{\begin{tabular}[t]{c}load angle $\theta_\mathrm{L}$ in rad\end{tabular}}}}}%
    \put(0,0){\includegraphics[width=\unitlength,page=3]{plots/ApplExRefTrack.pdf}}%
    \put(0.7425,0.52225){\makebox(0,0)[lt]{\lineheight{1.25}\smash{\begin{tabular}[t]{l}reference\end{tabular}}}}%
    \put(0,0){\includegraphics[width=\unitlength,page=4]{plots/ApplExRefTrack.pdf}}%
    \put(0.7425,0.4835){\makebox(0,0)[lt]{\lineheight{1.25}\smash{\begin{tabular}[t]{l}1st period\end{tabular}}}}%
    \put(0,0){\includegraphics[width=\unitlength,page=5]{plots/ApplExRefTrack.pdf}}%
    \put(0.7425,0.44475){\makebox(0,0)[lt]{\lineheight{1.25}\smash{\begin{tabular}[t]{l}5th period\end{tabular}}}}%
    \put(0,0){\includegraphics[width=\unitlength,page=6]{plots/ApplExRefTrack.pdf}}%
    \put(0.7425,0.406){\makebox(0,0)[lt]{\lineheight{1.25}\smash{\begin{tabular}[t]{l}10th per.\end{tabular}}}}%
    \put(0,0){\includegraphics[width=\unitlength,page=7]{plots/ApplExRefTrack.pdf}}%
    \put(0.7425,0.36725){\makebox(0,0)[lt]{\lineheight{1.25}\smash{\begin{tabular}[t]{l}50th per.\end{tabular}}}}%
    \put(0,0){\includegraphics[width=\unitlength,page=8]{plots/ApplExRefTrack.pdf}}%
    \put(0.7425,0.3285){\makebox(0,0)[lt]{\lineheight{1.25}\smash{\begin{tabular}[t]{l}100th per.\end{tabular}}}}%
    \put(0,0){\includegraphics[width=\unitlength,page=9]{plots/ApplExRefTrack.pdf}}%
    \put(0.7425,0.28975){\makebox(0,0)[lt]{\lineheight{1.25}\smash{\begin{tabular}[t]{l}optimal\end{tabular}}}}%
    \put(0,0){\includegraphics[width=\unitlength,page=10]{plots/ApplExRefTrack.pdf}}%
  \end{picture}%
\endgroup%

%% file: plots/ApplExLcurve.pdf_tex
\begingroup%
  \makeatletter%
  \providecommand\color[2][]{%
    \errmessage{(Inkscape) Color is used for the text in Inkscape, but the package 'color.sty' is not loaded}%
    \renewcommand\color[2][]{}%
  }%
  \providecommand\transparent[1]{%
    \errmessage{(Inkscape) Transparency is used (non-zero) for the text in Inkscape, but the package 'transparent.sty' is not loaded}%
    \renewcommand\transparent[1]{}%
  }%
  \providecommand\rotatebox[2]{#2}%
  \newcommand*\fsize{\dimexpr\f@size pt\relax}%
  \newcommand*\lineheight[1]{\fontsize{\fsize}{#1\fsize}\selectfont}%
  \ifx\svgwidth\undefined%
    \setlength{\unitlength}{195bp}%
    \ifx\svgscale\undefined%
      \relax%
    \else%
      \setlength{\unitlength}{\unitlength * \real{\svgscale}}%
    \fi%
  \else%
    \setlength{\unitlength}{\svgwidth}%
  \fi%
  \global\let\svgwidth\undefined%
  \global\let\svgscale\undefined%
  \makeatother%
  \begin{picture}(1,0.92307692)%
    \lineheight{1}%
    \setlength\tabcolsep{0pt}%
    \put(0,0){\includegraphics[width=\unitlength,page=1]{plots/ApplExLcurve.pdf}}%
    \put(0.2,0.091795){\makebox(0,0)[t]{\lineheight{1.25}\smash{\begin{tabular}[t]{c}0\end{tabular}}}}%
    \put(0.43461538,0.091795){\makebox(0,0)[t]{\lineheight{1.25}\smash{\begin{tabular}[t]{c}50\end{tabular}}}}%
    \put(0.66923077,0.091795){\makebox(0,0)[t]{\lineheight{1.25}\smash{\begin{tabular}[t]{c}100\end{tabular}}}}%
    \put(0.90384615,0.091795){\makebox(0,0)[t]{\lineheight{1.25}\smash{\begin{tabular}[t]{c}150\end{tabular}}}}%
    \put(0.55192346,0.02230769){\makebox(0,0)[t]{\lineheight{1.25}\smash{\begin{tabular}[t]{c}period\end{tabular}}}}%
    \put(0,0){\includegraphics[width=\unitlength,page=2]{plots/ApplExLcurve.pdf}}%
    \put(0.18564115,0.13076923){\makebox(0,0)[rt]{\lineheight{1.25}\smash{\begin{tabular}[t]{r}0.08\end{tabular}}}}%
    \put(0.18564115,0.24871808){\makebox(0,0)[rt]{\lineheight{1.25}\smash{\begin{tabular}[t]{r}0.1\end{tabular}}}}%
    \put(0.18564115,0.36666654){\makebox(0,0)[rt]{\lineheight{1.25}\smash{\begin{tabular}[t]{r}0.12\end{tabular}}}}%
    \put(0.18564115,0.48461538){\makebox(0,0)[rt]{\lineheight{1.25}\smash{\begin{tabular}[t]{r}0.14\end{tabular}}}}%
    \put(0.18564115,0.60256423){\makebox(0,0)[rt]{\lineheight{1.25}\smash{\begin{tabular}[t]{r}0.16\end{tabular}}}}%
    \put(0.18564115,0.72051269){\makebox(0,0)[rt]{\lineheight{1.25}\smash{\begin{tabular}[t]{r}0.18\end{tabular}}}}%
    \put(0.18564115,0.83846154){\makebox(0,0)[rt]{\lineheight{1.25}\smash{\begin{tabular}[t]{r}0.2\end{tabular}}}}%
    \put(0.05615385,0.50000038){\rotatebox{90}{\makebox(0,0)[t]{\lineheight{1.25}\smash{\begin{tabular}[t]{c}mean error per period\end{tabular}}}}}%
    \put(0,0){\includegraphics[width=\unitlength,page=3]{plots/ApplExLcurve.pdf}}%
    \put(0.51923077,0.68384615){\makebox(0,0)[lt]{\lineheight{1.25}\smash{\begin{tabular}[t]{l}with learning\end{tabular}}}}%
    \put(0,0){\includegraphics[width=\unitlength,page=4]{plots/ApplExLcurve.pdf}}%
    \put(0.51923077,0.62807692){\makebox(0,0)[lt]{\lineheight{1.25}\smash{\begin{tabular}[t]{l}without learning\end{tabular}}}}%
    \put(0,0){\includegraphics[width=\unitlength,page=5]{plots/ApplExLcurve.pdf}}%
    \put(0.51923077,0.57230769){\makebox(0,0)[lt]{\lineheight{1.25}\smash{\begin{tabular}[t]{l}optimal\end{tabular}}}}%
    \put(0,0){\includegraphics[width=\unitlength,page=6]{plots/ApplExLcurve.pdf}}%
  \end{picture}%
\endgroup%

%% file: plots/AppInputs.pdf_tex
\begingroup%
  \makeatletter%
  \providecommand\color[2][]{%
    \errmessage{(Inkscape) Color is used for the text in Inkscape, but the package 'color.sty' is not loaded}%
    \renewcommand\color[2][]{}%
  }%
  \providecommand\transparent[1]{%
    \errmessage{(Inkscape) Transparency is used (non-zero) for the text in Inkscape, but the package 'transparent.sty' is not loaded}%
    \renewcommand\transparent[1]{}%
  }%
  \providecommand\rotatebox[2]{#2}%
  \newcommand*\fsize{\dimexpr\f@size pt\relax}%
  \newcommand*\lineheight[1]{\fontsize{\fsize}{#1\fsize}\selectfont}%
  \ifx\svgwidth\undefined%
    \setlength{\unitlength}{195bp}%
    \ifx\svgscale\undefined%
      \relax%
    \else%
      \setlength{\unitlength}{\unitlength * \real{\svgscale}}%
    \fi%
  \else%
    \setlength{\unitlength}{\svgwidth}%
  \fi%
  \global\let\svgwidth\undefined%
  \global\let\svgscale\undefined%
  \makeatother%
  \begin{picture}(1,0.92307692)%
    \lineheight{1}%
    \setlength\tabcolsep{0pt}%
    \put(0,0){\includegraphics[width=\unitlength,page=1]{plots/AppInputs.pdf}}%
    \put(0.2,0.091795){\makebox(0,0)[t]{\lineheight{1.25}\smash{\begin{tabular}[t]{c}0\end{tabular}}}}%
    \put(0.55192308,0.091795){\makebox(0,0)[t]{\lineheight{1.25}\smash{\begin{tabular}[t]{c}5\end{tabular}}}}%
    \put(0.90384615,0.091795){\makebox(0,0)[t]{\lineheight{1.25}\smash{\begin{tabular}[t]{c}10\end{tabular}}}}%
    \put(0.55192346,0.02230769){\makebox(0,0)[t]{\lineheight{1.25}\smash{\begin{tabular}[t]{c}time within period in s\end{tabular}}}}%
    \put(0,0){\includegraphics[width=\unitlength,page=2]{plots/AppInputs.pdf}}%
    \put(0.18564115,0.13076923){\makebox(0,0)[rt]{\lineheight{1.25}\smash{\begin{tabular}[t]{r}-50\end{tabular}}}}%
    \put(0.18564115,0.24871808){\makebox(0,0)[rt]{\lineheight{1.25}\smash{\begin{tabular}[t]{r}0\end{tabular}}}}%
    \put(0.18564115,0.36666654){\makebox(0,0)[rt]{\lineheight{1.25}\smash{\begin{tabular}[t]{r}50\end{tabular}}}}%
    \put(0.18564115,0.48461538){\makebox(0,0)[rt]{\lineheight{1.25}\smash{\begin{tabular}[t]{r}100\end{tabular}}}}%
    \put(0.18564115,0.60256423){\makebox(0,0)[rt]{\lineheight{1.25}\smash{\begin{tabular}[t]{r}150\end{tabular}}}}%
    \put(0.18564115,0.72051269){\makebox(0,0)[rt]{\lineheight{1.25}\smash{\begin{tabular}[t]{r}200\end{tabular}}}}%
    \put(0.18564115,0.83846154){\makebox(0,0)[rt]{\lineheight{1.25}\smash{\begin{tabular}[t]{r}250\end{tabular}}}}%
    \put(0.05615385,0.50000038){\rotatebox{90}{\makebox(0,0)[t]{\lineheight{1.25}\smash{\begin{tabular}[t]{c}input $u$ in V\end{tabular}}}}}%
    \put(0,0){\includegraphics[width=\unitlength,page=3]{plots/AppInputs.pdf}}%
    \put(0.61153846,0.80692308){\makebox(0,0)[lt]{\lineheight{1.25}\smash{\begin{tabular}[t]{l}constr. $V_\mathrm{max}$\end{tabular}}}}%
    \put(0,0){\includegraphics[width=\unitlength,page=4]{plots/AppInputs.pdf}}%
    \put(0.61153846,0.75115385){\makebox(0,0)[lt]{\lineheight{1.25}\smash{\begin{tabular}[t]{l}1st period\end{tabular}}}}%
    \put(0,0){\includegraphics[width=\unitlength,page=5]{plots/AppInputs.pdf}}%
    \put(0.61153846,0.69538462){\makebox(0,0)[lt]{\lineheight{1.25}\smash{\begin{tabular}[t]{l}10th period\end{tabular}}}}%
    \put(0,0){\includegraphics[width=\unitlength,page=6]{plots/AppInputs.pdf}}%
    \put(0.61153846,0.63961538){\makebox(0,0)[lt]{\lineheight{1.25}\smash{\begin{tabular}[t]{l}50th period\end{tabular}}}}%
    \put(0,0){\includegraphics[width=\unitlength,page=7]{plots/AppInputs.pdf}}%
    \put(0.61153846,0.58384615){\makebox(0,0)[lt]{\lineheight{1.25}\smash{\begin{tabular}[t]{l}100th per.\end{tabular}}}}%
    \put(0,0){\includegraphics[width=\unitlength,page=8]{plots/AppInputs.pdf}}%
    \put(0.61153846,0.52807692){\makebox(0,0)[lt]{\lineheight{1.25}\smash{\begin{tabular}[t]{l}150th per.\end{tabular}}}}%
    \put(0,0){\includegraphics[width=\unitlength,page=9]{plots/AppInputs.pdf}}%
    \put(0.61153846,0.47230769){\makebox(0,0)[lt]{\lineheight{1.25}\smash{\begin{tabular}[t]{l}optimal\end{tabular}}}}%
    \put(0,0){\includegraphics[width=\unitlength,page=10]{plots/AppInputs.pdf}}%
  \end{picture}%
\endgroup%

%% file: 5_app.tex
\section{Application: Servomechanism}\label{sec:App}
In the two previous examples, we have not considered the actual time available for the computation.
In fact, the learning update may not always be computable within one sampling period $T_\mathrm{s}$, especially when considering systems with higher dimensions and fast dynamics. 
To circumvent this issue, we will outsource learning onto a server that communicates with the controller.
In particular we assume, that the convex hull update is executed in parallel to the controller and whenever it is completed, one of the latest data points since the last one added is added next.
Further, $i_\mathrm{T}(k)$ is not predefined but the optimization iteration is stopped as soon as the computation time is exhausted.
This does not affect our stability and convergence results as long as at least one optimization iteration is performed since Theorem \ref{thm:stability} and Theorem \ref{thm:convergence} apply for any sequence $i_\mathrm{T}(k) \in \N_+$.

In this section, we will implement the learning scheme for a reference tracking task on a simulation model of a servomechanism consisting of a DC-motor, a gear-box, an elastic shaft and a load. This system has already served as an example in \cite{BEMPORAD1998} from where we have inherited the linear model
\begin{flalign}\label{eq:applsys}
\dot x &= \begin{bmatrix}
0 & 1 & 0 & 0 \\
\frac{-k_\theta}{J_\mathrm{L}}& 
\frac{-\beta_\mathrm{L}}{J_\mathrm{L}}& 
\frac{k_\theta}{\rho J_\mathrm{L}}& 0 \\
0&0&0&1\\
\frac{k_\theta}{\rho J_\mathrm{M}}& 0 &
\frac{-k_\theta}{\rho^2 J_\mathrm{M}}& \frac{-\beta_\mathrm{M}R-K_\mathrm{T}^2}{J_\mathrm{M}R}
\end{bmatrix} x + \begin{bmatrix}
0\\0\\0\\ \frac{K_\mathrm{T}}{R J_\mathrm{M}}
\end{bmatrix}u&
\end{flalign}
and the parameters, which can be found in \cite{BEMPORAD1998}.
The state vector $x=[\theta_\mathrm{L}\ \dot \theta_{\mathrm L}\  \theta_\mathrm{M}\  \dot \theta_\mathrm{M}]^\top$ consists of the load angle $\theta_{\mathrm L}$ and the motor angle $\theta_{\mathrm M}$ as well as their time derivatives and the input $u$ corresponds to the DC voltage.
The reference tracking task of this system is also included in the collection of benchmark MPC problems given in \cite{KOUZOUPIS2015}. 
The model is discretized using zero-order hold on the input and the sampling time $T_\mathrm{s} =\unit[0.1]{s}$. 
Further the system has to satisfy the state and input constraints $\left|\left[
		k_\theta\ 0\  \sfrac{-k_\theta} \rho\ 0
	\right]x \right|\leq T_\mathrm{max},\ 
	|u|\leq V_\mathrm{max}$, where $ T_\mathrm{max}, V_\mathrm{max}$ can be found in \cite{BEMPORAD1998}.
To make the problem more difficult and to obtain a polytopic feasible set, we added the constraints $|x_1|\leq 2,\ |x_2| \leq 2,$ and $|x_2+x_4|\leq 40$ to the original problem.

The goal is that the angle of the load $\theta_\mathrm{L}$ follows the periodic step reference signal $r(k)$ with period length $200$.
To achieve this we see that $x_\mathrm{s}(r_\mathrm{s})=[1\ \ 0\ \ \rho\ \ 0]^\top r_\mathrm{s},\ u_\mathrm{s}=0$ is a steady state of \eqref{eq:applsys} that produces exactly the constant reference $r_\mathrm{s}$ for arbitrary $r_\mathrm{s}\in\R$. 
Therefore we will define the cost function with respect to $x-x_s(r(k))$ and $u-u_\mathrm{s}=u$ except for the part that incorporates the constraints. 
That is, \eqref{eq:l} and \eqref{eq:F} become
\begin{subequations}
\begin{align}
l(x,r,u)&=\norm{x-x_\mathrm{s}(r)}_Q^2+\norm{u}_R^2  +\varepsilon \hat B(x,u),\\
F(x,r)&=\norm{x-x_\mathrm{s}(r)}_P^2,
\end{align}
\end{subequations}
with $Q=\mathrm{diag}\,([10\ 0.1 \ 10\ 0.1]),\ R=0.01$, and $ \varepsilon=10^{-3}$.
Further, the relaxed logarithmic barrier function $\hat B$ and $P$ are chosen according to \cite{FELLER2017b}.
For the optimization update operator, we use the gradient descent backtracking line search as described in \cite{FELLER2017} with parameters $\rho=0.5$, $c_1=10^{-3}$ and $c_2=0.999$. 
The data set $\D(0)$ will be initialized with the extreme points of the feasible set as described in Remark \ref{rem:initCH}. 

The simulation results of the closed loop over $N_\mathrm{sim}=30000$ time steps ($150$ periods of the reference signal) are shown in Fig. \ref{fig:applEx}. 
In the left subplot, we can see that the reference tracking performance improves from period to period. 
This is also visible in the middle subplot, where the mean reference tracking error per period of the learning scheme decreases and converges to the one of the optimally controlled system. 
The original anytime MPC scheme without a spatial warm start and thus without learning does not improve over time and the mean reference tracking error stays constant. 
The scheme with learning is in the first period worse than without learning, because the time to compute the spatial warm start cannot be used for further optimization iterations.
In average there is a difference of $1.8$ iterations in between the original anytime MPC ($\approx 39.4$ iterations) and the scheme with learning ($\approx 37.6$).
However, already in the second period the scheme with learning achieves better tracking performance and improves much more over time.
After the $120$th period the controller has learned to exploit the full range of feasible control inputs as can be seen in the right subplot of Fig. \ref{fig:applEx}.

There is a computational aspect that we have not discussed yet.
The required memory capacity and the computation time for updating the convex hull and for the spatial warm start generation increase with the size of $\D(k)$ such that the number of points in the convex hull that can be handled is limited.
Hence, to exploit the available capacities we add only data points $(U,x,J)$ to $\D(k)$ that lead to a significant improvement from $J_N^\mathrm{a}(x,\D(k))$ to $J=J_N^\mathrm{a}(x,\D(k+1))$, here $> 10^{-2}$. 
In this example, out of $30000$ data points $21164$ were skipped because the improvement was too small and $4630$ because it was not possible to update the convex hull within one sampling time, leading to a total of $4206$ data points added to the convex hull.

In summary, the example clearly shows the benefit of the proposed memory-based MPC scheme in improving performance, while maintaining guarantees on stability at all times and being suitable for online implementation.

%% file: 6_con.tex
\section{Conclusion and outlook}\label{sec:Conc}
We presented and analyzed an online learning scheme for the value function and optimal policy in a real-time MPC framework for both linear and nonlinear systems. 
Even if only one optimization iteration is performed between two consecutive sampling instants, the MPC still approaches optimality in the long run through the learning scheme.
Furthermore, stability of the real-time MPC scheme is retained independent of the type of learning approach. 
These findings were illustrated in different linear and nonlinear numerical examples.
Moreover, we discussed how learning and applying the learned input can be decoupled and parallelized.

An interesting topic for future study is to further improve the presented learning methods by more efficient implementations of the data structure and the algorithms to generate the spatial warm start. 
In particular, while we have not provided a concrete algorithm for the nonlinear Lipschitz based learning scheme, the ideas of \cite{BELIAKOV2006} might be helpful for arriving at an efficient implementation. 
The optimization iteration also leaves room for improvement by exploiting underlying properties or structure. 
The proposed method can be extended in several ways.
For example the online learning scheme could probably be modified to handle time-varying systems or costs by including some mechanisms to forget expired data points.
Since our stability result does not depend on the learning method, this opens the possibility to apply classical machine learning function approximators for online learning.
A comparison to a function approximation based on offline data to generate a static spatial warm start would also be interesting to investigate.

%% file: A_proofs.tex
\section{Proofs}\label{sec:appendix_proofs}
\begin{pf}[Proof of Lemma \ref{lem:JNaProps}]
	(i) Since $J_N$ is convex and the spatial warm start defined in \eqref{eq:psi_sw_lin} is a convex combination of data points $(U_i,x_i,J_i)\in\mathcal D(k)$, we have
	\begin{align*}
	\big (\Psi_{\mathrm{sw}}(x,\mathcal D(k)),x,J_N^{\mathrm{a}}(x,\mathcal D(k))\big )&=\sum_{i=1}^K \lambda_i (U_i,x_i,J_i)&
	\end{align*}
	and further
	\begin{align*}\nonumber
	J_N\big (\Psi_\mathrm{sw}\big (x,\D(k)\big ),x\big ) &= J_N\left( \sum_{i=1}^K \lambda_i (U_i,x_i) \right)\\\nonumber &\leq  \sum_{i=1}^K \lambda_i J_N(U_i,x_i) \\&= \sum_{i=1}^K \lambda_iJ_i = J_N^\mathrm{a} \big (x,\D(k)\big ).
	\end{align*}
	
	(ii) We slightly reformulate \eqref{eq:psi_sw_lin} 
	\begin{align}\label{eq:JNaReformulated}
	J_N^\mathrm{a} \big (x,\D(k)\big )=\min_{\substack{(z,J)\in \conv \Dxj(k)\\ z=x}} J
	\end{align}
	and show that the epigraph of $J_N^\mathrm{a} (\cdot ,\D(k)):\conv \Dx(k)\to \R_+$ is convex:
	\begin{align*}
	\epi J_N^\mathrm{a}\big (\cdot, &\D(k)\big )= \left\{ (x,\bar J)\left| \bar J \geq\min _{\substack{(z, J)\in \conv \Dxj(k)\\ z= x}} J  \right\}\right.\\
	&= \left\{  (x,\bar J)\left|\exists J \leq\bar J: (x, J) \in \conv\Dxj(k) \right\}\right.\\
	&= \left\{  (x, J+\alpha) \left|\alpha \geq 0\wedge ( x, J) \in \conv\Dxj(k) \right\}\right.\\
	&= \conv \Dxj (k) \oplus \big ( \{0\} \times [0,\infty) \big )
	\end{align*}
	where $\oplus$ denotes the Minkowski sum. Since both $\conv \Dxj (k)$ and $ \{0\} \times [0,\infty)$ are convex, the epigraph $\epi J_N^\mathrm{a}(\cdot ,\D(k))$ is also.		
\end{pf}
\begin{pf}[Proof of Lemma \ref{lem:JNconv}]
It is shown in \cite{FELLER2017} that the cost function $J_N$ can be written under these assumptions as
	\begin{align*}
	J_N(U,x)=\frac 1 2 U^\top H U + x^\top F U + \frac 1 2 x^\top Y x + \varepsilon \hat B_\mathrm{xu} (U,x)
	\end{align*}
	where $H,F$ and $Y$ can be computed from $Q$, $R$, $P$, $A$ and $B$ and account for the terms $\norm{x}_Q^2$, $\norm{u}_R^2$ and $\norm{x}_P^2$ in $l$ and $F$. Furthermore, the relaxed logarithmic barrier function $\hat B_\mathrm{xu}(U,x)$ is shown in \cite{FELLER2017} to be a weighted sum over convex functions whose argument depends affinely on $(U,x)$. Therefore $\hat B_\mathrm{xu}(U,x)$ is convex in $(U,x)$. The quadratic part of $J_N$ is also convex in $(U,x)$ since it is quadratic and positive definite as proven in \cite{FELLER2017}.
\end{pf}
\begin{pf}[Proof of Lemma \ref{lem:JNaPropsNL}]
	In view of \eqref{eq:Lipsch} we have
	\begin{align*}
	&\big | J_N (\Psi_\mathrm{sw}(x,\D),x)- J_N \big (\Psi_\mathrm{sw}(x,\D),x^*(x,\D)\big ) \big |\\&\qquad\qquad\qquad\qquad \leq L(\Psi_{\mathrm{sw}}(x,\D)) \norm {x^*(x,\D) -x}
	\end{align*}
	and further
	\begin{align*} J_N \big (\Psi_\mathrm{sw}(x,\D),x\big )&\leq J_N \big (\Psi_\mathrm{sw}(x,\D),x^*(x,\D)\big ) \\&\qquad + L\big (\Psi_{\mathrm{sw}}(x,\D)\big ) \norm {x^*(x,\D)-x}  \\&= J_N^\mathrm{a} (x,D).
	\end{align*}
	To prove the second statement, we use \eqref{eq:xJstar} with $\bar x=y\in\R^n$ and \eqref{eq:JNagenau} with $\bar x=x\in\R^n$ to see
	\begin{align*}
	J_N^\mathrm{a} (x,\mathcal D)&\leq J^*(y,\D)+L(\Psi_\mathrm{sw}(y,\D)) \norm{x-x^*(y,\D)}
	\end{align*}
	which implies
	\begin{align*}
	&J_N^\mathrm{a} (x,\mathcal D)- J_N^\mathrm{a} (y,\mathcal D)\\ &\ \leq L\big (\Psi_\mathrm{sw}(y,\D)\big )\big ( \norm{x-x^*(y,\D)} - \norm{y-x^*(y,\D)}\big )\\ &\ \leq L_\D \norm{x-y}.
	\end{align*}
	Since $x$ and $y$ are interchangeable, it follows 
	\begin{align*}
	\big | J_N^\mathrm{a} (x,\mathcal D)- J_N^\mathrm{a} (y,\mathcal D) \big| & \leq L_\D \norm{x-y}. \ \square
	\end{align*}
\end{pf}
\begin{pf}[Step 2) of the proof of Theorem \ref{thm:convergence_nonl}]
	Step 2) $J_N^\mathrm{a}(x,\D(k))$ converges to $J_N^*(x)$: For $x\in \Omega$, there exists a sequence $k_i \to \infty$ such that $\xi_{i}=f(x(k_i),\Pi_0 U(k_i)) \to x$. It follows
	\begin{align*}
		\big|J_N^\mathrm{a} \big (\xi_i, \D(k_{i}&+1)\big )-J_N^\mathrm{a}\big (\xi_i,\D(k_i)\big ) \big| \leq 2M \norm{\xi_i-x}\\&+ \big|J_N^\mathrm{a} \big (x, \D(k_{i}+1)\big )-J_N^\mathrm{a}\big (x,\D(k_i)\big ) \big| \to 0.
	\end{align*}
		Thus the following chain of inequalities
		\begin{align}\label{eq:sandwich3}
		\begin{split}
		J_N^{\mathrm a} \big (\xi_i,&\D(k_i+1)\big )\leq J_N\big (U(k_i+1),\xi_i\big )\\&\leq J_N\big (\Phi^1\big (U(k_i),x(k_i),\D(k_i)\big ),\xi_i\big )\\&\leq J_N\big (\Phi^0\big (U(k_i),x(k_i),\D(k_i)\big ),\xi_i\big )\\&\leq J_N\big (\Psi_\text{sw}\big (\xi_i,\D(k_i)\big ),\xi_i\big )\\& \leq J_N^{\mathrm a}\big (\xi_i,\D(k_i)\big ),
		\end{split}
		\end{align}
		is a chain of equalities in the limit, where the first inequality holds due to \eqref{eq:D} and \eqref{eq:psi_sw_nonl} and the other inequalities due to \eqref{eq:ineqo}, \eqref{eq:ineqo}, \eqref{eq:psi_w} and Lemma \ref{lem:JNaPropsNL} in this specific order. Therefore as $i\to\infty$ the decrease of the optimizer update operator \eqref{eq:ineqo} 
		\begin{align*}
		\gamma \left( \Phi^0\big (U(k_i),x(k_i),\D(k_i)\big ), \xi_i \right) \to 0,
		\end{align*}
		which is only possible if 
		\begin{align*}
		J_N \left( \Phi^0\big (U(k_i),x(k_i),\D(k_i)\big ), \xi_i \right)-J_N^*(\xi_i) \to 0.
		\end{align*}
		This is due to \eqref{eq:sandwich3} equivalent to
		\begin{align*}
		J_N^\mathrm{a} \big ( \xi_i,\D(k_i) \big )-J_N^*(\xi_i) \to 0.
		\end{align*}
		$J_N^*$ is Lipschitz continuous and without loss of generality we can assume that $M$ is the Lipschitz constant of $J_N^*$ (if not choose $M$ large enough). Hence it follows
		\begin{align*}
		\big | J_N^\mathrm{a}& \big ( x,\D(k_i) \big ) -J_N^*(x) \big | \leq \big | J_N^\mathrm{a} \big ( \xi_i ,\D(k_i) \big ) -J_N^*(\xi_i) \big |\\
		&+\big | J_N^\mathrm{a} \big ( x ,\D(k_i) \big )-J_N^\mathrm{a} \big ( \xi_i ,\D(k_i) \big ) + J_N^*(\xi_i)-J_N^*(x) \big |\\ &\leq  \big | J_N^\mathrm{a} \big ( \xi_i ,\D(k_i) \big ) -J_N^*(\xi) \big |+2M \norm{x-\xi_i} \to 0
		\end{align*}
		which is with \eqref{eq:sandwich} what we wanted.
\end{pf}

%% file: A_algo.tex
\section{Convex hull algorithm}\label{A:algorithm}
\newcommand{\CH}{\mathcal{CH}}
\begin{table}[t]
	\caption{Properties of the convex hull object $\CH$}\label{tab:propCH}
	\begin{tabularx}{\linewidth}{cX}
		\toprule
		name & description \\ \midrule
		$\Dx$ & list of previous data points in $\R^n$\\
		$\Du$ & list of previous input sequences in $\R^{Nm}$\\
		$\Dj$ & list of previous cost function values in $\R_+$ \\
		$\CH_{\mathrm{xJ}}$ & list that contains all facets of the lower half of $\conv \Dxj$ as lists of the indices of their extreme points in $\Dxj$\\
		$\CH_{\mathrm{x}}$ & list that contains all facets of $\conv \Dx$ as lists of the indices of their extreme points in $\Dx$ \\
		$c_{\mathrm{x}}$ & center point of $\conv \Dx$ \\
		$c_{\mathrm{xJ}}$ & center point of $\conv \Dxj$ \\
		$o_{\mathrm{x}}$ & list of deleted facets in $\CH_{\mathrm{x}}$ that can be overwritten \\
		$o_{\mathrm{xJ}}$ & list of deleted facets in $\CH_{\mathrm{xJ}}$ that can be overwritten \\
		$\mathcal{G}_{\mathrm{xJ}}$ & list that contains for each point in $\Dxj$ a list of facets in $\CH_{\mathrm{xJ}}$ that are attached to this point \\
		$\mathcal{G}_{\mathrm{x}}$ & list that contains for each point in $\Dx$ a list of facets in $\CH_{\mathrm{x}}$ that are attached to this point \\
		\bottomrule
	\end{tabularx}
\end{table}
The purpose of the convex hull algorithm presented in this section is to solve \eqref{eq:psi_sw_lin}. As discussed in Section \ref{sec:linMPC}, $J_N^\mathrm{a}(\cdot,\D(k))$ gives rise to a triangulation of $\conv \Dx(k)$ and the spatial warm start at $x$ can be computed by seeking for the $n$-simplex in this triangulation that contains $x$. Therefore in a first step, this triangulation of $\conv \Dx(k)$ must be computed by a convex hull algorithm that determines the facets of the lower boundary of the convex hull $\conv \Dxj(k)$. The lower boundary of $\conv \Dxj(k)$ is the graph of $J_N^\mathrm{a} (\cdot,\D(k))$, see Fig. \ref{fig:psi_sw_lin}. The extreme points of each of these facets are, after projection onto $\R^n$, the extreme points of an $n$-simplex of the triangulation of $\conv \Dx(k)$. As second step, the spatial warm start can be computed from this convex hull by searching the $n$-simplex that contains $x$ and combine the inputs corresponding to its extreme points to obtain a warm start solution.
Hence, the procedure naturally decomposes into two steps:
\begin{enumerate}[label={\roman*)}]
	\item \textit{Generate spatial warm start} from the convex hull.
	\item \textit{Update convex hull} with a new data point.
\end{enumerate}
While i) corresponds to line 4 of Algorithm \ref{algo:anyLearn}, ii) corresponds to line 14 where collecting data is meant as fitting the new point into the data structure. Notice that this decomposition divides the algorithm into learning ii) and applying the learned spatial warm start i) -- a property that will be exploited in Section \ref{sec:App} for parallelization.

It is efficient to use an incremental convex hull algorithm that updates the existing convex hull object when a new data point arrives instead of starting the calculation from scratch. Further, the fact that for a point that is added to the convex hull, a spatial warm start was already generated and it has therefore been located in the triangulation can be leveraged for efficiency. The problem of searching for the $n$-simplex in the triangulation that contains $x$ is similar to the location problem in explicit MPC and hence these algorithms (see e.g. \cite{TONDEL2003}) can be used. The implementation we used, however, is based on a directed search by tracking neighboring $n$-simplices starting from an initial guess down to $x$. 

We present a way to implement the two main routines spatial warm start generation \texttt{generateSW} and convex hull update \texttt{updateCH}. First, we need to define an object that represents the collected data and all we need to know about the convex hull, we call this the convex hull object and denote it with $\mathcal{CH}$. The properties of the convex hull object are listed in Table \ref{tab:propCH}. References to a property of the convex hull object are denoted with a dot, e.g. $\CH.\Dx$. In addition to these listed properties we will use $\CH.\D$ and $\CH.\Dxj$ to denote the combined lists of the data points, so the $k$th element of $\CH.\D$ is for example a triple containing the $k$th elements of $\CH.\Du$, $\CH.\Dx$ and $\CH.\Dj$. 

The way the data is generated allows for efficient tailored implementations of the two routines, if the neighbors of each facet are known.
\begin{enumerate}[label={\roman*)}]
	\item The facet that contains $x$ can be found by starting from some initial facet and tracking neighboring facets towards $x$, for example by following $a+s(x-a)$ as $s$ increases from $0$ to $1$, where $a$ is some point in the facet. A good initial guess can be given by the facet in which the last point was located, if the system state does not change too fast from one time instant to the next. If this point was added to the convex hull, then this facet does not exist anymore. In this case any facet attached to this newly added point can be taken as initial guess. Another more advanced initial guess can be provided by saving for every point in the convex hull its subsequent point. If the subsequent point was not added, then an extreme point of the facet where the subsequent point was located is stored. This information provides for each point in the convex hull a good guess of its subsequent point and this can be used as basis for the initial guess. Thereby, the system dynamic gets also included implicitly in the prediction from the initial guess.
	\item A data point that needs to be added to the convex hull is already located in the triangulation of $\conv \Dx(k)$ by the previous warm start generation. Hence, by starting from this facet and checking the neighboring facets successively, all facets that need to be updated can be found.
\end{enumerate}
The information about neighboring facets are stored in the graph $\mathcal{CH}.\mathcal{G}_\mathrm{xJ}$, which needs to be updated whenever the convex hull is. 

\subsection{Generating the spatial warm start solution from the convex hull object}
Let $x$ be the point where we want to generate the warm start solution and $\CH$ the convex hull object. The main part in the warm start generation is the search for the facet that contains $x$. As discussed above, it is easy to provide a good guess for a point $a\in\mathcal{CH}.\Dx$ that is close to $x$. Assume the index $i$ of $a$ in $\mathcal{CH}.\Dx$ is given, then the routine takes the inputs $\CH$, $x$ and $i$. As outputs it makes sense to not only return the spatial warm start $\Psi_{\mathrm{sw}}=\Psi_{\mathrm{sw}}(x,\CH.\D)$ but also the index $F$ of the facet of $\CH.\CH_{\mathrm{xJ}}$ that contains $x$ and the index $i$ of the next initial guess. If $x$ is not contained in any facet, i.e. if $x\notin \conv \CH.\Dx$, then let $F$ be negative, in fact let it be $F=-F_\mathrm{x}$ where $F_\mathrm{x}$ is the index of a facet of $\CH.\CH_\mathrm{x}$ that excludes $x$ by means of $x$ lies outside the supporting halfspace to $\conv \Dx$ at this facet $F_\mathrm{x}$. Summing, the routine can be invoked by $(\Psi_{\mathrm{sw}},F,i)=\texttt{generateSW}(\CH,x,i)$. The basic idea in finding the facet that contains $x$ is to start the search at $a$ and track the facets along the line $a+s(x-a)$, $s\in[0,1]$ until $s=1$, then we have reached the facet that contains $x$. If we sometime step out of $\conv \mathcal{CH}.\Dx$, then $x\notin \conv \mathcal{CH}.\Dx$ and we cannot provide a spatial warm start as described in \eqref{eq:psi_sw_lin} at $x$, but we can at the point where we stepped out of  $\conv \mathcal{CH}.\Dx$ and use this as spatial warm start solution for $x$. The procedure can be described by the following steps, where the numbers indicate to which lines of pseudo-code in Algorithm \ref{algo:generateSW} the step corresponds to
\begin{enumerate}[label={\alph*)}]
	\item {\small (1:$-$6:)} Get the initial point $a$ and a list $f$ of indices of the facets in $\CH.\CH_{\mathrm{xJ}}$ that are attached to it.
	\item {\small (7:$-$16:)} Find the facet $F\in f$ in which the vector $v=x-a$ points starting from $a$.
	\item {\small (17:$-$27:)} If no such facet exists, then $v$ must point out of $\conv \CH.\Dx$, $a$ must lie on the boundary and hence $x\notin \conv \CH.\Dx$. So there must exist at least one facet of $\CH.\CH_{\mathrm{x}}$ attached to $a$ that excludes $x$. Find the index of this facet, return its negative as $F$ and return the input used at $a$ as spatial warm start. Also return the index of the extreme point of this facet that is closest to $x$ increased by one as next initial guess $i$. 
	\item {\small (28:$-$34:)} Else track the facets along $a+sv$ until $s\geq 1$ or no further facet exists along the line because we stepped out of $\conv \CH.\Dx$.
	\item {\small(35:$-$41:)} If we stepped out before $s=1$, then find the index of facet of $\CH.\CH_{\mathrm{x}}$ where we stepped out and return its negative as $F$, further return the index of the extreme point of this that is closest to $x$ increased by one as next initial guess $i$ and the convex combination of the inputs of this facet as spatial warm start $\Psi_{\mathrm{sw}}$.
	\item {\small (42:$-$46:)} Else we have found index $F$ of the facet of $\CH.\CH_{\mathrm{xJ}}$ that contains $x$ and return it. Also return the convex combination of the inputs at this facet as spatial warm start and the index of the extreme point of this that is closest to $x$ increased by one as next initial guess $i$.
\end{enumerate}
\begin{algorithm}
	\caption{Spatial warm start generation}\label{algo:generateSW}
	$(\Psi_{\mathrm{sw}},F,i)=$ \texttt{generateSW}($\mathcal{CH},x,i$)\\
	\textbf{Input:} convex hull object $\CH$, evaluation point $x$ and initial guess $i$\\
	\textbf{Output:} spatial warm start $\Psi_\mathrm{sw}$ at $x$, index $F$ of facet $\CH.\CH_{\mathrm{xJ}}$ that contains $x$ and next initial guess $i$. If $x\notin \conv \CH.\Dx$, then $F$ is negative and $|F|$ is the index of a facet in $\CH.\CH_\mathrm{x}$ that excludes $x$.
	\vspace{-0.2cm}
	\begin{enumerate}[label={\small \arabic*:}]\setlength{\itemsep}{-1ex}
		\item $f=\mathcal{CH}.\mathcal{G}_\mathrm{xJ}(i)$
		\item  \textbf{while} $f<0$ \textbf{do}
		\item \quad $i=-f$
		\item \quad $f=\mathcal{CH}.\mathcal G_\mathrm{xJ}(i)$
		\item \textbf{end while}
		\item $a=\mathcal{CH}.\Dx(i)$
		\item $v=x-a$
		\item $s=0$
		\item \textbf{for} $F\in f$ \textbf{do}
		\item \quad \textbf{if} \texttt{pointsInFacet}$(\CH,x,F,i)$ \textbf{then}
		\item \qquad $E=\CH.\CH_\mathrm{xJ}(F)$ without $i$
		\item \qquad $d=$ \texttt{normal}$(\CH.\Dx(E),a)$
		\item \qquad $s= d^\top(\CH.\Dx(E(1))-a)/(d^\top v)$
		\item \quad \quad \textbf{break for}
		\item \quad \textbf{end if}
		\item \textbf{end for}
		\item \textbf{if} $s=0$ \textbf{then}
		\item \quad \textbf{for} $F\in \CH.\mathcal G_\mathrm{x}(i)$
		\item \qquad $d= \texttt{normal}(\CH.\Dx(\CH.\CH_{\mathrm{x}}(F)),\CH.c_\mathrm{x})$
		\item \qquad \textbf{if} $d^\top (x-\CH.\Dx(\CH.\CH_{\mathrm{x}}(F)(1)))>0$ \textbf{then }
		\item \qquad \quad \textbf{break for} 
		\item \qquad \textbf{end if}
		\item \quad \textbf{end for}
		\item \quad $f=\CH.\CH_{\mathrm{x}}(F)$
		\item \quad $i=\argmin_{j\in f} \norm{x-\CH.\Dx(j)}$
		\item \quad \textbf{return} $(\Psi_{\mathrm{sw}},F,i)=(\CH.\Du(i),-F,i+1)$
		\item \textbf{end if}
		\item \textbf{while} true 
		\item \quad \textbf{if} $s\geq 1$ \textbf{then break while}
		\item \quad $f=$ \texttt{getFacets}$(\CH,E,\mathrm{xJ})$
		\item \quad $F=f$ without $F$
		\item \quad \textbf{if} $F$ is empty \textbf{then break while} 
		\item \quad $(s,E)=$ \texttt{findIntersection}$(\CH,a,v,s,F,E)$ 
		\item \textbf{end while}
		\item \textbf{if} $s<1$ \textbf{then}
		\item \quad $F=$ \texttt{getFacets}$(\CH,E,\mathrm{x})$
		\item \quad $i=\argmin_{j\in E} \norm{x-\CH.\Dx(j)}$
		\item \quad $c=$ \texttt{findConvComb}$(\CH,E,a+sv)$
		\item \quad $\Psi_{\mathrm{sw}}=\sum_{j=1}^n \CH.\Du(E(j)) c(j)$
		\item \quad \textbf{return} $(\Psi_{\mathrm{sw}},F,i)=(\Psi_{\mathrm{sw}},-F,i+1)$
		\item \textbf{end if}
		\item $f=\CH.\CH_{\mathrm{xJ}}(F)$
		\item $i=\argmin_{j\in f} \norm{x-\CH.\Dx(j)}$
		\item $c=$ \texttt{findConvComb}$(\CH,f,x)$
		\item $\Psi_{\mathrm{sw}}=\sum_{j=1}^{n+1} \CH.\Du(f(j)) c(j)$
		\item \textbf{return} $(\Psi_{\mathrm{sw}},F,i)=(\Psi_{\mathrm{sw}},F,i+1)$
	\end{enumerate}\vspace{-0.2cm}
\end{algorithm}
A detailed description of the steps in pseudo-code can be found in Algorithm \ref{algo:generateSW}, still there are needed some minor supporting routines that we will discuss now to complete the whole procedure of generating the warm start. 
\begin{itemize}
	\item  $b=\texttt{pointsInFacet}(\CH,x,F,i)$ checks if the vector $v=x-\CH.\Dx(i)$ starting from $i$ points inside the facet of $\CH.\CH_{\mathrm{xJ}}$ with index $F$, where $i$ has to be an extreme point of $F$. The length of $v$ does not matter only the direction, i.e. $x=\CH.\Dx(i)+v$ need not lie inside $F$ for $b$ being true, it is enough if there exists $\varepsilon>0$ such that $\CH.\Dx(i)+\varepsilon v$ lies inside the facet. Pseudo-code is given in Algorithm \ref{algo:pIF}. 
	\item $d=\texttt{normal}(A,x)$ is the most used subroutine and it calculates for a given set $A$ of $p$ points in $\R^p$ a vector $d\in\R^p$ that is normal to the hyperplane going through all points in $A$. This hyperplane cuts $\R^p$ into two halfspaces and $d$ points into the halfspace that does not contain the given point $x\in\R^p$, where $x$ must not lie on the this plane. A pseudo-code description of this subroutine is given in Algorithm \ref{algo:normal} and uses the well-known QR-decomposition that decomposes an invertible matrix $M\in\R^{p\times p}=QR$ into an orthogonal matrix $Q=Q^{-\top}$ and an upper triangular matrix $R$ where all diagonal entries of $R$ are positive. 
	\item $f=\texttt{getFacets}(\CH,E,\mathrm{flag})$ gives out all facets of $\CH.\CH_{\mathrm{flag}}$ that share the extreme points specified in $E$. $f$ can also be only a single facet or even empty, depending on how many facets exists that share the edge $E$. Pseudo-code is given in Algorithm \ref{algo:gf}. 
	\item $f=\texttt{findIntersection}(\CH,a,v,s,F,E)$ takes an edge $E$ of $F$ such that $a+sv$ lies on $E$ and it gives out $(s,E)$ such that $a+sv$ lies on the edge $E$ of $F$ but with output $s>$ input $s$. In words it takes the edge and point where the line $a+sv$ enters $F$ and computes the edge and point where it leaves $F$. Pseudo-code is given in Algorithm \ref{algo:fI}. 
	\item $f=\texttt{fincConvComb}(\CH,f,x)$ calculates the weights $c$ such that $\sum_{j=1}^{|f|} c(j) \CH.\Dx(f(j))=x$ and $\sum_{j=1}^{|f|} c(j)=1$. Therefore it must be verified that it is possible to convex combine the points in $f$ to $x$. Pseudo-code is given in Algorithm \ref{algo:fcc}. 
\end{itemize}
\begin{algorithm}
\caption{Check if vector points in facet}\label{algo:pIF}
$b=$ \texttt{pointsInFacet}($\mathcal{CH},x,F,i$)\\
\textbf{Input:} convex hull object $\CH$, vector $x$, facet $F$ and index $i$ of starting point \\
\textbf{Output:} boolean $b$, true if $v=x-\CH.\Dx(i)$ points in $F$ starting from $i$.\vspace{-0.2cm}
\begin{enumerate}[label={\small \arabic*:}]\setlength{\itemsep}{-1ex}
	\item $E=\CH.\CH_\mathrm{xJ}(F)$
	\item \textbf{for} $j=1,2,\dots,n+1$
	\item \quad \textbf{if} $E(j)\neq i$ \textbf{then}
	\item \qquad $\bar E=E\text{ without }E(j)$
	\item \qquad $d=$ \texttt{normal}$(\CH.\Dx(\bar E),\CH.\Dx(E(j)))$
	\item \qquad \textbf{if} $d^\top (x-\CH.\Dx(i))>0$ \textbf{then}
	\item \qquad\quad \textbf{return} $b=$ false
	\item \qquad \textbf{end if}
	\item \quad \textbf{end if}
	\item \textbf{end for}
	\item \textbf{return} $b=$ true
\end{enumerate}\vspace{-0.2cm}
\end{algorithm}

\begin{algorithm}
\caption{Calculate normal vector}\label{algo:normal}
$d=$ \texttt{normal}$(A,x)$\\
\textbf{Input:} list $A$ of $n$ points in $\R^n$, point $x\in\R^n$\\
\textbf{Output:} vector $d$ that is normal to the plane spanned by $A$ pointing into the halfspace that does not contain $x$\vspace{-0.2cm}
\begin{enumerate}[label={\small \arabic*:}]\setlength{\itemsep}{-1ex}
\item \textbf{for} $j=1,\dots,n-1$
\item \quad $A(j)=A(j)-A(n)$
\item \textbf{end for}
\item $A(n)=x-A(n)$
\item QR-decomposition $QR=[A(1), A(2),\dots,A(n)]$
\item \textbf{return} $d=-$ last column of $Q$
\end{enumerate}\vspace{-0.2cm}
\end{algorithm}

\begin{algorithm}
	\caption{Get facets attached to set of extreme points}\label{algo:gf}
	$f=$ \texttt{getFacets}$(\CH,E,\mathrm{flag})$\\
	\textbf{Input:} convex hull object $\CH$, list $E$ of extreme points and $\mathrm{flag}$ that can either be $\mathrm{x}$ or $\mathrm{xJ}$ indicating the convex hull \\
	\textbf{Output:} list $f$ of facets of $\CH.\CH_{\mathrm{flag}}$ that contain all extreme points given in $E$\vspace{-0.2cm}
	\begin{enumerate}[label={\small \arabic*:}]\setlength{\itemsep}{-1ex}
		\item $f=\CH.\mathcal{G}_\mathrm{flag}(E(1))$
		\item \textbf{for} $p\in E$ without $E(1)$ \textbf{do}
		\item \quad $f=$ common elements of $f$ and $\CH.\mathcal{G}_\mathrm{flag}(p)$
		\item \textbf{end for}
		\item \textbf{return} $f$
	\end{enumerate}\vspace{-0.2cm}
\end{algorithm}

\begin{algorithm}
\caption{Find intersection of line and facet boundary}\label{algo:fI}
$(s,E)=$ \texttt{findIntersection}$(\mathcal{CH},a,v,s,F,E)$\\
\textbf{Input:} convex hull object $\CH$, start point $a$, direction vector $v$, current $s$, facet $F$, edge $E$ of $F$ with $a+sv$ lying on this edge \\
\textbf{Output:} $s$ and $E$ such that $a+sv$ lies on ende $E$ of facet $F$ and output $s>$ input $s$, \vspace{-0.2cm}
\begin{enumerate}[label={\small \arabic*:}]\setlength{\itemsep}{-1ex}
	\item $p=\CH.\CH_\mathrm{xJ}(F)$ without $E$
	\item \textbf{for} $j=1,2,\dots,n$
	\item \quad $\bar E(j)=E\text{ without }E(j)\text{ but with }p$
	\item \quad $d=$ \texttt{normal}$(\CH.\Dx(\bar E),\CH.\Dx(E(j)))$
	\item \quad $\bar s(j)=d^\top(\CH.\Dx(p)-a)/(d^\top v)$
	\item \quad \textbf{if} $\bar s(j)\leq s$ \textbf{then} $\bar s(j)=\infty$
	\item \textbf{end for}
	\item $\bar j=\argmin_j \bar s(j)$ 
	\item \textbf{return} $(s,E)=(\bar s(\bar j), \bar E(\bar j))$
\end{enumerate}\vspace{-0.2cm}
\end{algorithm}

\begin{algorithm}
	\caption{Find convex combination}\label{algo:fcc}
	$c=$ \texttt{findConvComb}$(\CH,E,x)$\\
	\textbf{Input:} convex hull object $\CH$, index array $E$ of points in $\CH.\Dx$, point $x$\\
	\textbf{Output:} vector $c$ such that the points in $E$ summed up with the weights in $c$ equal $x$ and sum over elements in $c$ is $1$ \vspace{-0.2cm}
	\begin{enumerate}[label={\small \arabic*:}]\setlength{\itemsep}{-1ex}
		\item $A=[\CH.\Dx(E)^\top\ \mathbbm{1}^\top]^\top$
		\item $b=[x^\top \ 1]^\top$
		\item Solve $Ac=b$ for $c$
		\item \textbf{return} $c$
	\end{enumerate}\vspace{-0.2cm}
\end{algorithm}

\subsection{Updating the convex hull}

Second we consider updating the convex hull, therefore let $\CH$ be the convex hull object and $(U,x)$ the input sequence and point in state space we want to add. In addition we already know the index $F$ of the facet of $\CH.\CH_{\mathrm{xJ}}$ that contains $x$ from the previous spatial warm start generation at exactly this point $x$. This makes up the inputs of the routine and the output is of course the updated convex hull object $\CH$, hence it can be invoked by $\CH=\texttt{updateCH}(\CH,x,U,F)$. The basic idea of updating $\CH$ is that we start at the facet $F$, which must be for sure updated since there is always an improvement in the optimization update operator as long as we have not reached exact optimality and even in that unlikely case the point we want to add lies directly on the facet and it makes no difference to add it and split the facet. Thus we always update the facet $F$ and the we start searching from $F$ for neighboring facets that need to be updated and if we found one, then we also search in its neighbors for facets to update. Since the convex hull was convex before and has to be convex afterwards the set of facets to update must be connected and thus we will find all facets to update by that procedure. A special case is if $x\notin \conv \CH.\Dx$, i.e. $F<0$, then we need to find all facets of $\CH.\CH_{\mathrm{x}}$ that exclude $x$ and check that facets of $\CH_{\mathrm{xJ}}$ attached to them if they must be updated. In this case we must also update $\CH.\CH_{\mathrm{x}}$ which we do by the same procedure, starting from $-F$ search for neighboring facets that must to be updated. In fact all facets of $\CH.\CH_{\mathrm{x}}$ that must be updated are facets that exclude $x$ and the other way round, so we search them, update them and initialize with the facets of $\CH.\CH_{\mathrm{xJ}}$ attached to them the search for updates in $\CH.\CH_{\mathrm{xJ}}$. In more detail the procedure can be described by the following steps, where the numbers indicate to which lines of pseudo-code in Algorithm \ref{algo:update} the step corresponds to
\begin{algorithm}
	\caption{Update convex hull object}\label{algo:update}
	$\CH=$ \texttt{updateCH}$(\CH,x,U,F)$\\
	\textbf{Input:} convex hull object $\CH$, new point $x$, input sequence $U$, index $F$ of facet in $\CH.\CH_\mathrm{xJ}$ that contains $x$, negative if $x\notin \conv \CH.\Dx$, then $|F|$ is facet of $\CH.\CH_\mathrm{x}$ that excludes $x$.\\
	\textbf{Output:} updated convex hull object $\CH$ that contains the new point\vspace{-0.2cm}
	\begin{enumerate}[label={\small \arabic*:}]\setlength{\itemsep}{-1ex}
		\item Append $x$ to $\CH.\Dx$
		\item Append $U$ to $\CH.\Du$
		\item Append $J_N(U,x)$ to $\CH.\Dj$
		\item $z=[x^\top\ J_N(U,x)]^\top$
		\item $k=$ number of elements in $\Dx$
		\item \textbf{if} $F<0$ \textbf{then}
		\item \quad $\mathcal E_\text{toTestx}=$ list of all edges of facet $\CH.\CH_{\mathrm{x}}(-F)$
		\item \quad $\mathcal E_\text{toTestxJ}=$ list with one element $\CH.\CH_{\mathrm{x}}(-F)$
		\item \quad $\CH=\texttt{removeFacet}(\CH,-F,\mathrm{x})$
		\item \quad \textbf{while} $\mathcal E_\text{toTestx}$ is not empty \textbf{do}
		\item \qquad $F=$ \texttt{getFacets}$(\CH,\mathcal E_\text{toTestx}(1),\mathrm{x})$
		\item \qquad $E=\mathcal E_\text{toTestx}(1)$
		\item \qquad \textbf{if} $F$ is not empty \textbf{then}
		\item \qquad \quad $d=\texttt{normal}(\CH.\Dx(\CH.\CH_{\mathrm{x}}(F)),\CH.c_\mathrm{x})$
		\item \qquad \quad \textbf{if} $d^\top (x-\CH.\Dx(E(1)))>0$ \textbf{then}
		\item \qquad\qquad append $\CH.\CH_{\mathrm{x}}(F)$ to $\mathcal E_\text{toTestxJ}$
		\item \qquad\qquad $\CH=\texttt{removeFacet}(\CH,F,\mathrm{x})$
		\item \qquad\qquad append edges of $F$ except $E$ to $\mathcal E_\text{toTestx}$
		\item \qquad \quad \textbf{else}
		\item \qquad \qquad $\CH=\texttt{addFacet}(\CH,[E,\ k],\mathrm{x})$
		\item \qquad \quad \textbf{end if}
		\item \qquad \textbf{end if}
		\item \qquad remove $E$ from $\mathcal E_\text{toTestx}$
		\item \quad \textbf{end while}
		\item \textbf{else}
		\item \quad $\mathcal E_\text{toTestxJ}=$ list of all edges of facet $\CH.\CH_{\mathrm{xJ}}(F)$
		\item \quad $\CH=\texttt{removeFacet}(\CH,F,\mathrm{xJ})$
		\item \textbf{end if}
		\item \textbf{while} $\mathcal E_\text{toTestxJ}$ is not empty \textbf{do}
		\item \quad $E=\mathcal E_\text{toTestxJ}(1)$
		\item \quad $F=$ \texttt{getFacets}$(\CH,E,\mathrm{xJ})$
		\item \quad \textbf{if} $F$ is not empty \textbf{then}
		\item \qquad $d=\texttt{normal}(\CH.\Dxj(\CH.\CH_{\mathrm{xJ}}(F)),\CH.c_\mathrm{xJ})$ 
		\item \qquad \textbf{if} $d^\top (z-\CH.\Dxj(E(1)))>0$ \textbf{then}
		\item \qquad\quad $\CH=\texttt{removeFacet}(\CH,F,\mathrm{xJ})$
		\item \qquad\quad append edges of $F$ except $E$ to $\mathcal E_\text{toTestxJ}$
		\item \qquad \textbf{else}
		\item \qquad \quad $\CH=\texttt{addFacet}(\CH,[E,\ k],\mathrm{x})$
		\item \qquad \textbf{end if}
		\item \quad \textbf{end if}
		\item \quad remove $E$ from $\mathcal E_\text{toTestxJ}$
		\item \textbf{end while}
		\item $\CH.c_\mathrm{x}=\CH.c_\mathrm{x} (k-1)/k+x/k$
		\item $\CH.c_\mathrm{xJ}=\CH.c_\mathrm{xJ} (k-1)/k+z/k$
		\item \textbf{return} $\CH$
	\end{enumerate}\vspace{-0.2cm}
\end{algorithm}

\begin{enumerate}[label={\alph*)}]
	\item {\small (1:$-$5:)} Append data point to data set $\CH.\D$.
	\item {\small (6:$-$24:)} If $F<0$ then $x\notin \conv \CH.\Dx$ and we need to update $\CH.\CH_{\mathrm{x}}$ first. We know that $-F$ is one facet that needs to be updated so we remove it and check all neighboring facets if they must also be updated. If so, then we remove them and also add their edges to the list of edges whose neighboring facets need to be checked. If not then we take the edge between the facet that had been removed and the one that stays and add it together with the new point $k$ as new facet to $\CH.\CH_\mathrm{x}$. Whenever we remove a facet we also add it to the list of edges we need to check for updating $\CH.\CH_{\mathrm{xJ}}$.
	\item {\small (25:$-$28)} If $F>0$ then we remove this facet from $\CH.\CH_{\mathrm{xJ}}$ and initialize the list of edges we need to check with all edges of this facet.
	\item {\small (29:$-$42)} Check edges as long as there are some whether their adjacent facets need to be updated or not, if so remove them and add their edges to the list, if not add the edge with the new point as new facet to $\CH.\CH_{\mathrm{xJ}}$.
	\item {\small (43:$-$45:)} Update the center points of the convex hulls and return the updated convex hull object.
\end{enumerate}
Whenever we talked about adding a facet to or removing it from $\CH.\CH_{\mathrm{x/xJ}}$ we of course also have to update the graph $\CH.\mathcal{G}_{\mathrm{x/xJ}}$. A detailed explanation on the subroutines adding and removing shall be given now.
\begin{itemize}
	\item $\CH=\texttt{removeFacet}(\CH,F,k,\mathrm{flag})$ removes the facet $F$ from the $\CH.\CH_{\mathrm{flag}}$ and also removes its appearance in $\CH.\mathcal{G}_\mathrm{flag}$. When removing $F$ we still need to keep the slot $F$ in $\CH.\CH_{\mathrm{flag}}$ because if we would delete it, all facets coming after $F$ in the list $\CH.\CH_{\mathrm{flag}}$ would get their index decremented and we would need to change the whole $\CH.\mathcal G_\mathrm{flag}$. Hence it is easier to add $F$ to a list $\CH.o_\mathrm{flag}$ of open slots that can be overwritten, thereby keeping all other indices as they are. If it happens for $\mathrm{flag}=\mathrm{xJ}$ that thereby an extreme point $p$ of $F$ is afterwards not attached to any facet at all, i.e. with adding the new data point $p$ 'moves' from the boundary of the convex hull to its interior, then we store for $p$ the information that it has fallen out of the boundary when updating $k$ by setting $\CH.\mathcal{G}_\mathrm{xJ}(p)=-k$. This is necessary for the initial guess for the warm start generation where it could happen that for some point the initial guess is $p$, but if $p$ is not attached to any facet then we cannot start a search from there. Therefore we store that $p$ has been overwritten by $k$ and then we can start the search from a facet attached to $k$ instead.
	\item $\CH=\texttt{addFacet}(\CH,f,\mathrm{flag})$ adds the facet with extreme points listed in $f$ to $\CH.\CH_{\mathrm{flag}}$ and updates the graph $\CH.\CH_{\mathrm{flag}}$. If there is an open slot $\CH.o_{\mathrm{flag}}$ in $\CH.\CH_{\mathrm{flag}}$ then we can overwrite it, otherwise we append the new facet to the list $\CH.\CH_{\mathrm{flag}}$.
\end{itemize}
Another point we have not discussed here is the initialization of the convex hull. The algorithms presented here only work if there already exists a convex hull object where there are enough points to build the convex hulls and they have to be correct obviously. A straight forward approach is to wait until the $n+1$ data points are available and initialize the convex hull $\CH.\CH_{\mathrm{xJ}}$ with the first and only facet $[1\ \dots\ n+1]$ as well as the convex hull $\CH.\CH_{\mathrm{x}}$ with all $n+1$ possible combinations of $n$ points out of these $n+1$ points. 

\begin{algorithm}
	\caption{Remove facet from convex hull object}\label{algo:remove}
	$\CH=$ \texttt{removeFacet}$(\CH,F,k,\mathrm{flag})$\\
	\textbf{Input:} convex hull object $\CH$, facet $F$ of $\CH.\CH_\mathrm{flag}$, current index $k$ and $\mathrm{flag}$ that is either $\mathrm{x}$ or $\mathrm{xJ}$ indicating the convex hull \\
	\textbf{Output:} updated convex hull object $\CH$ without $F$\vspace{-0.2cm}
	\begin{enumerate}[label={\small \arabic*:}]\setlength{\itemsep}{-1ex}
		\item \textbf{for} $p\in \CH.\CH_{\mathrm{flag}}(F)$ \textbf{do}
		\item \quad remove $F$ from $\CH.\mathcal{G}_\mathrm{flag}(p)$
		\item \quad \textbf{if} $\CH.\mathcal{G}_\mathrm{flag}(p)$ is empty and $\mathrm{flag}=\mathrm{xJ}$ \textbf{then}
		\item \qquad $\CH.\mathcal{G}_\mathrm{flag}(p)=-k$
		\item \quad \textbf{end if}
		\item \textbf{end for}
		\item add $F$ to $\CH.o_\mathrm{flag}$
		\item \textbf{return} $\CH$
	\end{enumerate}\vspace{-0.2cm}
\end{algorithm}

\begin{algorithm}
\caption{Add facet to convex hull object}\label{algo:add}
$\CH=$ \texttt{addFacet}$(\CH,f,\mathrm{flag})$\\
\textbf{Input:} convex hull object $\CH$, list $f$ of extreme points of facet that is added to $\CH.\CH_\mathrm{flag}$, current index $k$ and $\mathrm{flag}$ that is either $\mathrm{x}$ or $\mathrm{xJ}$ indicating the convex hull \\
\textbf{Output:} updated convex hull object $\CH$ with facet $f$\vspace{-0.2cm}
\begin{enumerate}[label={\small \arabic*:}]\setlength{\itemsep}{-1ex}
	\item \textbf{if} $\CH.o_\mathrm{flag}$ is not empty \textbf{then}
	\item \quad $F=\CH.o_\mathrm{flag}(1)$
	\item \quad change $\CH.\CH_\mathrm{flag}(F)$ to $f$
	\item \textbf{else}
	\item \quad $F=$ length of $\CH.\CH_\mathrm{flag}$ $+1$
	\item \textbf{end if}
	\item \textbf{for} $p\in f$ \textbf{do}
	\item \quad append $F$ to $\CH.\mathcal{G}(p)$
	\item \textbf{end for}
	\item remove $F$ from $\CH.o_\mathrm{flag}$
	\item \textbf{return} $\CH$
\end{enumerate}\vspace{-0.2cm}
\end{algorithm}

%% file: A_lipsch.tex
\section{Computation of the Lipschitz constant from Section \ref{sec:nonlAcEx}}\label{A:lipsch}
The Lipschitz constant $L(U)$ of the cost function $J_N(U,\cdot)$ can be computed as
\begin{align}
L(U)=\sum_{i=0}^{N}\big (L_l+\delta_{Ni}(L_F-L_l) \big )\prod_{j=0}^{i-1} L_f (\Pi_j U)
\end{align}
where $\delta_{Ni}$ denotes the Kronecker-delta that is $1$ if $i=N$ and $0$ else, where $L_l$, $L_F$ and $L_f$ are the Lipschitz constants of the stage cost $l$, the terminal cost $F$ and the system dynamic $f$, respectively and where $\Pi_j=[\dots\ 0, \ I_m,\ 0,\ \dots]\in \R^{m\times Nm}$ is a projection matrix, that projects $U$ onto its $mj+1$st till $mj+m$th component. The Lipschitz constants can be upper estimated as
\begin{subequations}
	\begin{align}
	L_l&=\sqrt{\frac {13}{50}},\qquad L_F=100 L_l\\ L_f(u)&=\sqrt{1+ T_\mathrm{s} |u_1| \sqrt{1+T_\mathrm{s}^2u_1^2/4} + T_\mathrm{s}^2u_1^2/2}.
	\end{align}	
\end{subequations}
by the following calculations 
\begin{align*}
\left(\deldel{l}{x_3}(x,u) \right)^2 &=\left( 0.1 \sin\left(\frac{x_3} 2\right)\cos\left(\frac{x_3} 2\right)\right)^2\leq 0.1^2\\
\left(\deldel{l}{x_{1,2}}(x,u) \right)^2 &=\left(\frac{ 2x_{1,2}}{4\left(1+x_1^2+x_2^2\right)^{\sfrac{3}{4}}}\right)^2 \\ &=  \frac{1}{4} \sqrt{\frac{x_{1,2}^4}{(1+x_1^2+x_2^2)^3}} \\&\leq  \frac{1}{4} \sqrt{\frac{x_{1,2}^4}{(1+x_{{1,2}}^2)^3}}\\
&=\frac{1}{4} \sqrt{\frac{x_{1,2}^4}{1+3 x_{1,2}^2+3x_{1,2}^4+x_{1,2}^6}}\leq \frac{1}{8}\\
\norm{ \deldel{l}{x}(x,u)}&\leq \sqrt{0.1^2+\frac{1}{4}}= \sqrt{\frac {13}{50}}
\end{align*}
and 
\begin{align*}
\deldel{f}{x}(x,u)&=\begin{bmatrix}
1 & 0 & -T_\mathrm{s} u_1\sin (x_3) \\0 & 1 & T_\mathrm{s} u_1 \cos(x_3) \\ 0&0&1
\end{bmatrix}\\
\Rightarrow\ \ \deldel{f}{x}(x,u)^\top \deldel{f}{x}(x,u)&=\\&\hspace{-2.5cm}\begin{bmatrix}
1 & 0 & -T_\mathrm{s} u_1\sin (x_3) \\0 & 1 & T_\mathrm{s} u_1 \cos(x_3) \\  -T_\mathrm{s} u_1\sin (x_3) & T_\mathrm{s} u_1 \cos(x_3) &T_\mathrm{s}^2 u_1^2 +1
\end{bmatrix}
\end{align*}
which has the eigenvalues $1$ and $1\pm T_\mathrm{s} u_1 \sqrt{1+T_\mathrm{s}^2u_1^2/4} + T_\mathrm{s}^2u_1^2/2$. Thus the maximum absolute eigenvalue is 
\begin{align*}
\norm{\deldel{f}{x}(x,u)}^2=1+ T_\mathrm{s} |u_1| \sqrt{1+T_\mathrm{s}^2u_1^2/4} + T_\mathrm{s}^2u_1^2/2.
\end{align*}